\documentclass[11pt, reqno]{amsart}
\usepackage{a4wide,color}
\usepackage{amsmath,amsthm,amsfonts,amssymb,mathtools}
\usepackage[colorlinks, citecolor = blue, draft = false]{hyperref}
\usepackage{tikz}
\usetikzlibrary{backgrounds}
\usetikzlibrary{math}
\usepackage{pgfplots}
\usetikzlibrary{fillbetween}
\usepackage{graphicx,tikz,subfigure,color,calc}
\usepackage{rotating}
\usepackage{crossreftools}
\usetikzlibrary{decorations.pathreplacing}
\usetikzlibrary{patterns}

\definecolor{darkgreen}{rgb}{0.0,0.5,0.0}
\definecolor{darkblue}{rgb}{0.0,0.0,1}
\definecolor{light-gray}{gray}{0.6}
\definecolor{light-red}{HTML}{b32c4b}
\definecolor{really-light-gray}{gray}{0.8}

\hypersetup{linkcolor = light-red}

\makeatletter
\renewcommand*{\eqref}[1]{%
  \hyperref[{#1}]{\textup{\tagform@{\ref*{#1}}}}%
}
\makeatother

\def\R{\mathbb R}

\def\N{\mathbb N}
\def\Z{\mathbb Z}
\def\P{\mathcal P}

\let\e\varepsilon
\DeclareMathOperator\dist{dist}

\DeclareMathOperator\supp{supp}
\DeclareMathOperator\curl{Curl}

\DeclareMathOperator\conv{conv}

\DeclareMathOperator{\SO}{SO}
\DeclareMathOperator{\Id}{Id}

\newtheorem{theorem}{Theorem}

\newtheoremstyle{remark}
{3pt}
{3pt}
{}
{}
{\bfseries}
{.}
{.5em}
{}
\theoremstyle{remark}
\newtheorem{remark}[theorem]{Remark}

\newcommand{\black}{\color{black}}

\title[Read--Shockley formula for a general Bravais lattice in two dimensions]{Read--Shockley formula for a general \\Bravais lattice in two dimensions}

\author[L.~Scardia]{L.~Scardia}
\author[E.G.~Tolotti]{E.G.~Tolotti}

\address[L.~Scardia]{Department of Mathematics, Heriot--Watt University, United Kingdom}\email{L.Scardia@hw.ac.uk} 

\address[E.G.~Tolotti]{Applied Mathematics M\"unster, University of M\"unster, Germany}\email{edoardo.tolotti@uni-muenster.de}

\begin{document}

 \begin{abstract}
        In this note we consider a two-dimensional semi-discrete dislocation energy and propose a simple and physically motivated construction for the grain boundary between two crystal grains with a small orientation difference. In the case of a general Bravais lattice, the energy of this construction matches the logarithmic scaling predicted by Read and Shockley.

                \bigskip


        \medskip

    \end{abstract}

    \maketitle

\section{Introduction and setting of the problem}
In this note we consider the two-dimensional nonlinear semi-discrete dislocation energy introduced by Lauteri and Luckhaus in \cite{LL} (and later studied and generalised in \cite{FGS}) and propose an alternative, simpler, and physically sound construction for the grain boundary between two crystal grains with a small orientation difference.
The energy of our construction is bounded by
\begin{equation*}
    \varepsilon E_{0}\theta(A - \log\theta),
\end{equation*}
 in agreement with the Read--Shockley formula in \cite{RS}, where $\varepsilon>0$ is a small parameter representing the microscopic scale (e.g.,  approximately the interatomic distance), and $\theta$ is the angle of misorientation. The constants $E_{0}$ and $A$ depend on the orientation of the grain boundary and on the lattice spacing; additionally, $A$ depends on the dislocation core radius, namely the size of the region around the defect where the use of a continuum model is not justified. 
We compute $E_{0}$ and $A$ explicitly, showing that they match the ones derived in \cite{RS} for a square lattice.

\subsection{The energy}  
Let $L>0$, and let $\Omega=[-L,L]\times [-2L,0]$ represent the two-dimensional cross-section of a three-dimensional crystal.
Let  $0 < \tau < \lambda$ be such that $\tau \varepsilon$ is the crystal lattice spacing and $\lambda \varepsilon$ is the dislocation core radius. 
Let $\mathcal{B}$ denote a Bravais lattice in $\mathbb{R}^2$ generated by $b_1, b_2 \in \mathbb{S}^{1}$, so that $\mathcal{B} = \text{span}_{\mathbb{Z}^2}\{b_1,b_2\}$, and $\text{span}_{\mathbb{R}^2}\{b_1,b_2\} = \R^2$. 
Finally, we denote with $\theta>0$ the small misorientation angle between two neighbouring grains occupying the regions $[-L,0]\times [-2L,0]$ and $[0,L]\times [-2L,0]$. 

The energy is defined on the admissible class 
$$
\mathcal{C}_\e(\Omega):=\left\{
(\beta,S): \beta \in L^1(\Omega; \R^{2\times 2}), S \subset \Omega \text{ relatively closed, satisfying \ref{regularity}--\ref{boundary-conditions}} 
\right\},
$$
where 
\begin{itemize}
\item[{\crtcrossreflabel{(H1)}[regularity]}] $\beta \in L^{2}(\Omega \setminus \overline{B_{\lambda \varepsilon}(S)})$ and $\supp (\curl \beta)\subset S$;
\item[{\crtcrossreflabel{(H2)}[circulation]}] for every simple, closed, and Lipschitz curve $\gamma \subset \Omega \setminus \overline{B_{\lambda \varepsilon}(S)}$ we have
    \begin{equation}\label{quant-circulation}
        \displaystyle \int_\gamma \beta t \,  d\mathcal{H}^{1} \in \tau \e \mathcal{B},
    \end{equation}
    where $t$ is the tangent vector\footnote{The orientation of the tangent vector $t$ is chosen taking the outer normal vector and rotating it counterclockwise by $\pi / 2$.} to $\gamma$;
\item[{\crtcrossreflabel{(H3)}[boundary-conditions]}] $\beta$ satisfies the symmetric boundary conditions
\begin{equation*}
\beta=
\begin{cases}
R_{-\theta} \quad &\text{in } [-L,-L+l]\times [-2L,0], \\
R_{\theta} \quad &\text{in } [L-l,L]\times [-2L,0], 
\end{cases}
\end{equation*}
with $0<l \ll L$, where $R_{\theta}$ is the matrix representing a counterclockwise rotation of angle $\theta$\footnote{Note that in \cite{RS} the boundary conditions are given in terms of $R_{\pm\theta/2}$. For notational simplicity we prefer to use $\theta$. We will comment on the differences arising from this choice whenever relevant.}, given by
\begin{equation*}
    R_{\theta} = 
    \begin{pmatrix}
        \cos\theta & -\sin\theta\\
        \sin\theta & \cos\theta
    \end{pmatrix}.
\end{equation*}
\end{itemize}
The rigorous interpretation of the integral in \eqref{quant-circulation} deserves some clarification.
By \ref{regularity}, we have that $\curl (\beta) = 0$ on $ \Omega \setminus \overline{B_{\lambda \varepsilon}(S)}$.
In particular, by classical arguments (see, for example \cite[Chapter IX, Thm. 2]{RJ}), one can prove that $\beta$ admits a tangent trace on every simple, closed, and Lipschitz curve $\gamma \subset \Omega \setminus \overline{B_{\lambda \varepsilon}(S)}$.
Then, \eqref{quant-circulation} can be interpreted in a weak sense as
\begin{equation*}
    \langle T^{t}_{\gamma}(\beta), 1 \rangle \in \tau \varepsilon \mathcal{B},
\end{equation*}
where $T_{\gamma}^{t}(\beta)$ is the tangent trace of $\beta$ along $\gamma$.

The energy associated to an admissible pair $(\beta,S)\in \mathcal{C}_\e(\Omega)$ is given by 
\begin{equation}\label{BL:energy}
    \mathcal{F}_\e(\beta,S) \coloneq \int_{\Omega\setminus B_{\lambda \e}(S)}\dist^2(\beta,\SO(2))dx + 
    \frac{\tau^{2}\black}{\lambda^2}\mathcal{L}^2(B_{\lambda \e}(S)) ,
\end{equation}
with the first term being the elastic energy of the body and the second term a core energy. The choice of the scaling of the core energy is for dimensional consistency, since the elastic energy carries a factor $\tau^2$, due to \eqref{quant-circulation}, which is compensated by the factor $\tau^{2}/\lambda^{2}$ in the core energy.

\subsection{Literature review and motivation}
The use of (nonlinear or linearised) semi-discrete energies in the mathematical modelling of dislocations is well established. The advantage of this mixed description, as opposed to a fully discrete or fully continuum formulation, is that it 
keeps track of each individual dislocation, while blurring the atomistic structure of the lattice. This is justified when the typical distance between dislocations is much larger than the atomistic spacing, which is the case for reasonable strains. 

The energy \eqref{BL:energy} is similar to the one introduced in \cite{MSZ-Indiana,SZ-SIMA}, which is in its turn the nonlinear version of the dislocation energy studied in \cite{GLP}. 
While in \cite{GLP,MSZ-Indiana,SZ-SIMA} it was assumed that dislocations were well separated, and hence their total number was bounded, here one needs the additional core energy term in \eqref{BL:energy} to ensure the same bound. A core energy term was needed also in e.g., in \cite{DLGP, Ginster}, where the separation condition was lifted, and in \cite{Ponsi}, in a slightly different framework. As observed there, this additional term is generally of smaller order compared with the elastic energy, and does not affect the leading-order asymptotics of the energy.

Semi-discrete energies have been recently used to describe small-angle tilt grain boundaries, and to predict the formation of polycrystals, see \cite{FPP, FGS, LL}. Small-angle tilt grain boundaries are well understood, and their energy is given by the celebrated Read--Shockley formula in \cite{RS}. In \cite{LL}, the validity of the formula was shown for the square lattice, starting from the energy \eqref{BL:energy}. Namely it was shown that, for $\mathcal{B}=\Z^2$,
\begin{equation}\label{aim:LL}
\inf\left\{
\mathcal{F}_\e(\beta,S): (\beta,S)\in \mathcal{C}_{\e}(\Omega)
\right\}
\leq C \varepsilon\theta (|\log \theta|+1),
\end{equation}
with $C>0$. The bound \eqref{aim:LL} was attained by means of an explicit construction consisting of an array of dislocations at the grain boundary, evenly spaced at distance $\varepsilon/\theta$. This result was later generalised in \cite{FGS} for any Bravais  lattice---albeit the explicit construction of the grain boundary was done only in the case of a possibly rotated square lattice.
Note that, allowing for a rotated lattice covers also the case of asymmetric boundary conditions for the square lattice $\mathcal{B} = \Z^{2}$, as considered in \cite{RS}. 
Indeed, in the latter case, it is sufficient to rotate the crystal to symmetrise the boundary conditions. 
The upper bound \eqref{aim:LL} for the energy was complemented by a matching lower bound in \cite{LL}, and by a compactness result for a special class of strains with bounded energy. In \cite{FGS}, a full compactness result was obtained, as well as lower and upper bounds for the energy in terms of $\Gamma$-convergence. In this note we focus on the grain boundary construction, and hence on providing the upper bound \eqref{aim:LL}. 

A linearised version of the Read--Shockley formula was derived in \cite{FPP} starting from a linearised version of the energy \eqref{BL:energy}, under the assumption of good separation for the dislocations, and for a diluted energy regime. Since the separation condition only allows for a number of dislocations much smaller than $\theta/\varepsilon$, dislocations cannot form a grain boundary as in \cite{RS}. The linearised energy in \cite{FPP} can then be used to describe `linearised polycrystals', where the misorientation of the grains is not only small, but infinitesimal.

Finally, the validity of the Read--Shockley formula has been recently investigated in \cite{GT22}, starting from a discrete model \`a la Ariza-Ortiz, and in \cite{GT25}, starting from an energy describing pairwise interactions of dislocations via a Volterra-like potential. Note that the energy in \cite{GT25} is still semi-discrete as \eqref{BL:energy}, namely it models individual dislocations in an otherwise continuum medium, even if it is a discrete sum.

\subsection{Our contribution}
In this work we deal with the case of symmetric boundary conditions for a general Bravais lattice $\mathcal{B}$, and propose a simple and physically motivated construction for the grain boundary, which we now describe briefly. We recall that the choice of symmetric boundary conditions is not restrictive since we can always reduce to this case by rotating the lattice.

In our construction, the admissible strain field $\beta$ agrees with the respective boundary conditions $R_{\pm \theta}$ in the majority of the domain, except for a thin vertical strip with width of order $\varepsilon/\theta$, as in \cite{FGS,LL}. 
In the general case where $e_1\notin \{b_1,b_2\}$,  
this strip is composed of two vertical sub-strips next to each other, and in each strip dislocations are evenly spaced, at a vertical distance of order  $\varepsilon / \theta$. 
The presence of a dislocation is incorporated by means of a jump of the deformation---and hence a nonzero circulation for the strain---that simulates the opening of the lattice to make space for the extra half plane of atoms carried by the defect.
We take into account the accumulation of the Burgers vector by making such opening wider and wider as we travel down in the lattice. Additionally, the circulation of the strain around each defect is along one generator of the lattice, $b_1$ or $b_2$, and only one of them is active in each sub-strip. This is a stronger condition than what is required in \ref{circulation}.  In the special case where  $e_1\in \{b_1,b_2\}$ only one strip is required, and the circulation of the strain is along $e_1$ whenever nonzero.  

We now highlight how our construction differs from the ones in \cite{FGS,LL}.
First of all, condition \ref{circulation} was first considered in this form in \cite{FGS}.
Its role is to impose that the circulation of the field $\beta$ is along a lattice direction, and that it is  quantized, in terms of the lattice spacing.
Note that, if $\beta$ is sufficiently regular, by Stokes' Theorem we have 
\begin{equation*}
    \int_{\gamma} \beta t \, d\mathcal{H}^{1} = \int_{\Gamma} \curl\beta \, dx,
\end{equation*}
where $\Gamma$ is the region enclosed by $\gamma$.
In particular, \ref{circulation} imposes that the averaged macroscopic Burgers vector lies in the lattice $\tau \varepsilon \mathcal{B}$.
In \cite{LL} only a weaker quantisation condition was required for the strain.

Moreover, we deal with any Bravais lattice, while in \cite{LL} only the square lattice with symmetric boundary conditions is considered.
In \cite{FGS} the explicit upper bound construction is done only in the case of a possibly rotated square lattice, although the lower bound works for any Bravais lattice.
Compared with the construction in \cite{FGS}, ours is simpler.
Indeed, the construction in \cite{FGS} is made up of three vertical sub-strips, rather than two. Moreover, our deformation is piecewise affine, and hence the strain is piecewise constant, and forms a polygonal Caccioppoli's partition of the domain, which is not the case in \cite{FGS}.

Our construction is in fact admissible for the more restrictive minimisation problem 
\begin{equation}\label{pc-beta}
    \inf\left\{
    \mathcal{F}_\e(\beta,S): (\beta,S)\in \mathcal{C}_{\e}(\Omega), \, \beta \text{ piecewise constant}
    \right\},
\end{equation}
which can be recast purely in terms of matrices.
Within this more restrictive framework, a good competitor for \eqref{pc-beta} is a piecewise constant map $\beta$ which connects the rotations $R_{\pm\theta}$ on the two outer vertical strips of the domain, in such a way that the constant values of $\beta$ are close to $\SO(2)$---to make the first term of the energy small---and rank-one connected along their common straight interfaces whenever possible---to make the second term of the energy small.
Note that, while the first term of the energy penalises strains that are far from rotations rather than not exact rotations, the second term penalises strains that are not rank-one connected.
Hence, to have a small energy, $\beta$ can afford not being in $\SO(2)$ in large portions of the domain (provided it is close to $\SO(2)$), while its values on neighbouring regions have to be rank-one connected on most of the common interfaces.
This forces the construction to have small, concentrated regions where the rank-one connection is not achieved, which can be identified with the locations of dislocations. 

We believe that this simpler construction is amenable to generalisations, e.g., to considering a multi-well energy and constraining the strain to be in $\SO(2)$ or some $\SO(2)U_i$, describing twins and laminates. These extensions are beyond the scope of this paper. 


\section{The grain boundary construction}

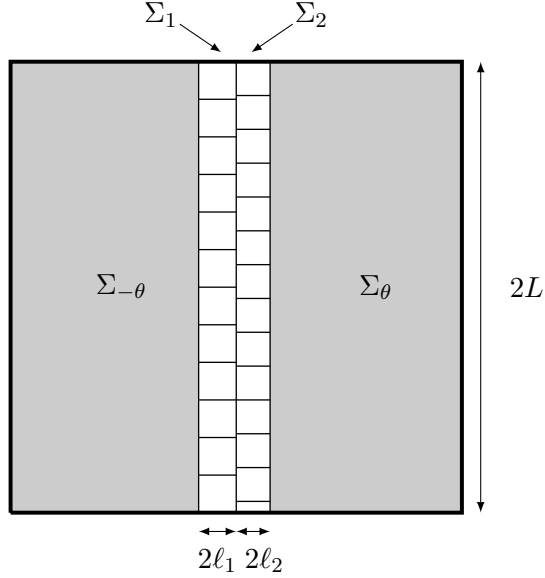
\begin{figure}[ht] 
	\begin{center}
		\begin{tikzpicture}[>=latex, scale=.25]

\fill[color=really-light-gray] (-12,-12)--(-2,-12)--(-2,12)--(-12,12);
\fill[color=really-light-gray] (12,-12)--(1.8,-12)--(1.8,12)--(12,12);

\draw[line width = 1.5pt](-12,-12)--(12,-12)--(12,12)--(-12,12)--(-12,-12);
\draw[<->] (13,-12) to (13,12);
\draw(14,0)node[right]{$2L$};
\draw[<->] (0,-13) to (1.8,-13);
\draw(1.5,-14.5)node[]{$2 \ell_2$};

\draw(-4,14.5)node[]{$\Sigma_1$};
\draw[->] (-3,14) to (-0.8,12.5);
\draw[->] (3,14) to (0.8,12.5);
\draw(4,14.5)node[]{$\Sigma_2$};
\draw[<->] (0,-13) to (-2,-13);
\draw(-1,-14.5)node[]{$2 \ell_1$};

\draw(-8,0)node[right]{$\Sigma_{-\theta}$};
\draw(6,0)node[right]{$\Sigma_{\theta}$};

\draw[line width = 0.5pt](0,-12)--(0,12);
\draw[line width = 0.5pt](1.8,-12)--(1.8,12);
\draw[line width = 0.5pt](0,-11.4)--(1.8,-11.4);
\draw[line width = 0.5pt](0,-9.6)--(1.8,-9.6);
\draw[line width = 0.5pt](0,-7.8)--(1.8,-7.8);
\draw[line width = 0.5pt](0,-6)--(1.8,-6);
\draw[line width = 0.5pt](0,-4.2)--(1.8,-4.2);
\draw[line width = 0.5pt](0,-2.4)--(1.8,-2.4);
\draw[line width = 0.5pt](0,-0.6)--(1.8,-0.6);
\draw[line width = 0.5pt](0,1.2)--(1.8,1.2);
\draw[line width = 0.5pt](0,3)--(1.8,3);
\draw[line width = 0.5pt](0,4.8)--(1.8,4.8);
\draw[line width = 0.5pt](0,6.6)--(1.8,6.6);
\draw[line width = 0.5pt](0,8.4)--(1.8,8.4);
\draw[line width = 0.5pt](0,10.2)--(1.8,10.2);

\draw[line width = 0.5pt](-2,-12)--(-2,12);
\draw[line width = 0.5pt](0,-10)--(-2,-10);
\draw[line width = 0.5pt](0,-8)--(-2,-8);
\draw[line width = 0.5pt](0,-6)--(-2,-6);
\draw[line width = 0.5pt](0,-4)--(-2,-4);
\draw[line width = 0.5pt](0,-2)--(-2,-2);
\draw[line width = 0.5pt](0,0)--(-2,0);
\draw[line width = 0.5pt](0,2)--(-2,2);
\draw[line width = 0.5pt](0,4)--(-2,4);
\draw[line width = 0.5pt](0,6)--(-2,6);
\draw[line width = 0.5pt](0,8)--(-2,8);
\draw[line width = 0.5pt](0,10)--(-2,10);
\end{tikzpicture}
	\end{center}
\caption{The splitting of the domain.}	\label{construction-Bravais}
\end{figure}

In this section we present our construction, illustrated in Figures \ref{construction-Bravais} and \ref{Fig:squares}.

\subsection{The splitting of the domain} 

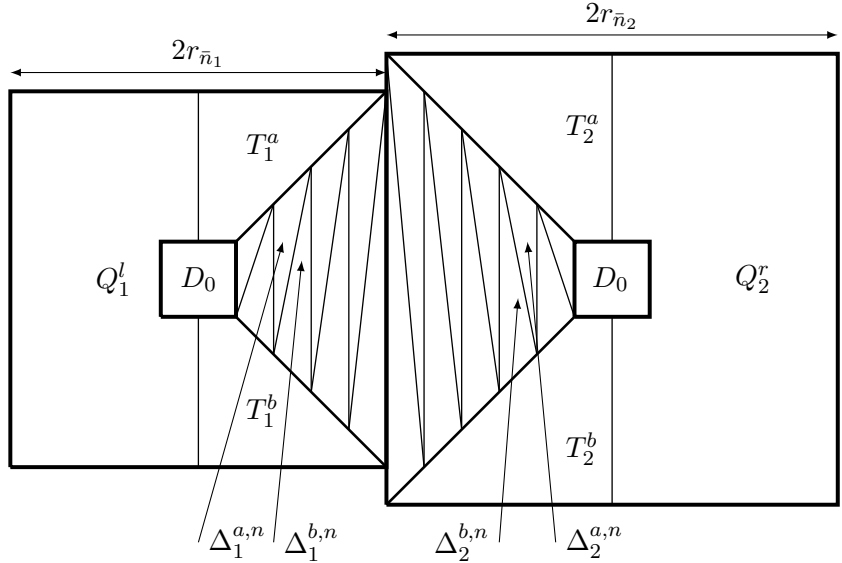
\begin{figure}[h] 
	\begin{center}
		\begin{tikzpicture}[>=latex, scale=.25]

\draw(10.4,0)node[right]{$D_{0}$};
\draw(18,0)node[right]{$Q_2^{r}$};
\draw(-11.5,0)node[right]{$D_{0}$};
\draw(-16,0)node[right]{$Q_1^l$};
\draw(9,8)node[right]{$T_2^{a}$};
\draw(9,-9)node[right]{$T_2^{b}$};

\draw(-8,7)node[right]{$T_1^{a}$};
\draw(-8,-7)node[right]{$T_1^{b}$};

\draw[<-] (-5.5,2) to (-10,-14);
\draw(-10,-14)node[right]{$\Delta_1^{a,n}$};

\draw[<-] (-4.5,1) to (-6,-14);
\draw(-6,-14)node[right]{$\Delta_1^{b,n}$};

\draw[<-] (7.5,2) to (9,-14);
\draw(9,-14)node[right]{$\Delta_2^{a,n}$};

\draw[<-] (7,-1) to (6,-14);
\draw(2,-14)node[right]{$\Delta_2^{b,n}$};

\draw[line width = 1.5pt](0,-12)--(24,-12)--(24,12)--(0,12)--(0,-12);
\draw[line width = 1.5pt](10,-2)--(14,-2)--(14,2)--(10,2)--(10,-2);

\draw[<->] (-20,11) to (0,11);
\draw(-10,12)node{$2 r_{\bar n_{1}}$};
\draw[line width = 1.5pt](-20,-10)--(0,-10)--(0,10)--(-20,10)--(-20,-10);
\draw[line width = 1.5pt](-12,-2)--(-8,-2)--(-8,2)--(-12,2)--(-12,-2);
\draw[line width = 1pt](-8,-2)--(0,-10);
\draw[line width = 1pt](-8,2)--(0,10);
\draw[line width = .5pt](-10,2)--(-10,10);
\draw[line width = .5pt](-10,-2)--(-10,-10);

\draw[line width = .5pt](-6,-4)--(-6,4);
\draw[line width = .5pt](-4,-6)--(-4,6);
\draw[line width = .5pt](-2,-8)--(-2,8);
\draw[line width = .5pt](-8,-2)--(-6,4);
\draw[line width = .5pt](-6,-4)--(-4,6);
\draw[line width = .5pt](-4,-6)--(-2,8);
\draw[line width = .5pt](-2,-8)--(0,10);

\draw[line width = .5pt](8,4)--(10,-2);
\draw[line width = .5pt](6,6)--(8,-4);
\draw[line width = .5pt](4,8)--(6,-6);
\draw[line width = .5pt](2,10)--(4,-8);
\draw[line width = .5pt](0,12)--(2,-10);

\draw[line width = .5pt](8,-4)--(8,4);
\draw[line width = .5pt](6,-6)--(6,6);
\draw[line width = .5pt](4,-8)--(4,8);
\draw[line width = .5pt](2,-10)--(2,10);

\draw[line width = 1pt](10,-2)--(0,-12);
\draw[line width = 1pt](10,2)--(0,12);

\draw[line width = .5pt](12,2)--(12,12);
\draw[line width = .5pt](12,-2)--(12,-12);

\draw[<->] (0,13) to (24,13);
\draw(12,14)node{$2 r_{\bar n_{2}}$};
\end{tikzpicture}
	\end{center}
\caption{The regions of the construction.}
\label{Fig:squares}
\end{figure}

As a first step, we split $\Omega$ into four vertical strips. 
Let $l<\ell_1, \ell_2\ll L$ be positive parameters to be defined later. We set 
$$
\Sigma_1 \coloneq [-2\ell_1,0]\times  [-2L,0], \quad \Sigma_2 \coloneq [0,2\ell_2]\times  [-2L,0],
$$
and 
$$
\Sigma_{-\theta} \coloneq [-L,-2\ell_1]\times [-2L,0], \quad \Sigma_\theta \coloneq [2\ell_2, L]\times [-2L,0],
$$
and we subdivide the domain as
\begin{align*}
\Omega=\Sigma_{-\theta}\cup \Sigma_1\cup \Sigma_2\cup \Sigma_\theta.
\end{align*}
Then, for $i=1,2$, we further split the strips $\Sigma_i$ into squares with side-length $2\ell_i$, obtained by translating $Q_i \coloneq [-\ell_i,\ell_i]^2$ 
by the vectors
\begin{equation}\label{translation-vector}
    t_{i,k} \coloneq -(2k-1)\ell_i e_2 + (-1)^i\ell_i e_1,
\end{equation}
where $k=1, \dots, N_i$, and $N_i \coloneq \lceil L/\ell_i \rceil$. We further partition $Q_i$ into a family of dyadic square annuli. To this aim, with no loss of generality we write $\ell_{i}$ as 
\begin{equation}\label{form-ell}
    \ell_{i} = 2^{\bar n_{i}} \lambda\varepsilon,
\end{equation}
where $\bar n_{i} \in \N$ will be fixed later, and define 
\begin{gather}
r_{0} \coloneq \lambda\varepsilon, \qquad r_{n} \coloneq 2^n r_{0},\label{squaresides}\\
D_{n}\coloneq 
\begin{cases}
  \medskip
    [-r_{0}, r_{0}]^{2}, &\quad n=0,\\ 
    \medskip
    [-r_{n},r_{n}]^2 \setminus [-r_{n-1},r_{n-1}]^2, &\quad n = 1, \dots, \bar n_{i}.
    \end{cases} \nonumber
\end{gather}
Then
$$
Q_i = [-\ell_{i}, \ell_{i}]^2 = [-r_{\bar n_{i}},r_{\bar n_{i}}]^2=\bigcup_{n=0}^{\bar n_i} D_{n}.
$$
We define the following subsets of $D_{n}$
\begin{align*}
\Delta_{i}^{a,n} & \coloneq \conv\{((-1)^{i+1}r_{n-1},-r_{n-1}), ((-1)^{i+1}r_{n},r_{n}),((-1)^{i+1}r_{n-1},r_{n-1})\},\\
\Delta_{i}^{b,n} & \coloneq \conv\{((-1)^{i+1}r_{n},r_{n}),((-1)^{i+1}r_{n},-r_{n}),((-1)^{i+1}r_{n-1},-r_{n-1})\},
\end{align*}
where $\conv$ denotes the convex hull.
Finally, we partition $Q_1$ and $Q_2$ as 
\begin{align}\label{partition1}
Q_1&=D_{0}\cup Q_1^l\cup T_1^{a}\cup T_1^{b}\cup \bigcup_{n=1}^{\bar n_{1}}\big(\Delta_1^{a,n}\cup \Delta_1^{b,n}\big),\\
\label{partition2}
Q_2&=D_{0}\cup Q_2^r\cup T_2^{a}\cup T_2^{b}\cup\bigcup_{n=1}^{\bar n_{2}} \big(\Delta_2^{a,n}\cup \Delta_2^{b,n}\big),
\end{align}
where $Q_1^l \coloneq [-r_{\bar n_{1}},0]^2\setminus D_{0}$, $Q_2^r \coloneq [0,r_{\bar n_{2}}]^2\setminus D_{0}$, and
\begin{align*}
T_i^{a}& \coloneq \conv\{(0,r_{0}),((-1)^{i+1}r_{0},r_{0}), ((-1)^{i+1}r_{\bar n_{i}},r_{\bar n_{i}}), (0, r_{\bar n_{i}})\},\\
T_i^{b}& \coloneq \conv\{(0,-r_{0}),((-1)^{i+1}r_{0}, -r_{0}), ((-1)^{i+1}r_{\bar n_{i}},-r_{\bar n_{i}}), (0, -r_{\bar n_{i}})\}.
\end{align*}
We refer to Figure \ref{Fig:squares} for an illustration of the partition of $Q_i$. Note that the sets of the partitions are convex polygons. 

The translated squares $Q^k_i:=Q_i+t_{i,k}$ in the strip $\Sigma_i$ are divided in the same way, and we indicate with a further subscript $k$ their partitioning sets, which are obtained by translating the sets in the right-hand side of \eqref{partition1}--\eqref{partition2} by the vector $t_{i,k}$ defined in \eqref{translation-vector}.

\subsection{The construction of the piecewise constant strain}\label{strain:con} In this section we construct a piecewise affine deformation in $\Omega$.
Note that the energy depends only on the strain field $\beta$, but we believe that writing the full deformation makes it easier to understand the construction.

In the vertical strips $\Sigma_{\pm
\theta}$ we define the deformation as $R_{\pm\theta}x$, for $x\in \Sigma_{\pm\theta}$, to match the boundary conditions \ref{boundary-conditions} for the strain. 

We now focus on the construction in the strips $\Sigma_i$. We will define a piecewise affine deformation in each of the squares $Q_i^k$, with constant gradients in the (translates of the) components of the partitions \eqref{partition1}--\eqref{partition2}. In a nutshell, we combine a rotation $R_i:=R_{\theta_i}$, with $\theta_i=(-1)^i\theta$, which is needed to match the boundary conditions on $\Sigma_{\pm\theta}$, with an additional translation, which is due to the presence of defects. More precisely, in the top regions $T_{i,k}^a$ we impose a translation along $b_i$ of $k$ units, in order to accommodate all the dislocations with Burgers vector $b_i$ up until that point.
In the bottom regions $T_{i,k}^b$, instead, we increment the translation of the top by one unit, to describe an additional dislocation. The lateral zig-zag regions are used to transition between the Burgers vector $b_1$ on $\Sigma_1$ and $b_2$ on $\Sigma_2$ in a way that does not create curl on their common boundary.

\smallskip

We use the following notation. 
For $p_i\in \R^2$ and $v_i\in \R^2$, $i=1,2,3$, we define the affine interpolation map $I_{\Delta}$ with values $v_i$ at the points $p_i$ as 
\begin{equation}\label{interpol}
    I_{\Delta}(x):= \sum_{i=1}^3 v_i \Phi_i(x), \qquad x\in \Delta:=\conv\{p_1,p_2,p_3\},
\end{equation}
where
\begin{align*}
& \Phi_i(x) \coloneq \Phi_i^{\text{ref}}(P^{-1}(x-p_1)), \quad P \coloneq (p_2-p_1, p_3-p_1) \in \R^{2 \times 2},\\
& \Phi_1^{\text{ref}}(x) \coloneq 1-x_1-x_2, \quad \Phi_2^{\text{ref}}(x) \coloneq x_1, \quad \Phi_3^{\text{ref}}(x) \coloneq x_2.
\end{align*}
The constant gradient of the interpolation map in $\Delta$ is 
\begin{equation}\label{const-grad}
    \begin{aligned}
    \nabla I_{\Delta}(x) &=  \sum_{i=1}^3 v_i \otimes \nabla \Phi_i(x)\\
    & = 
    v_1 \otimes P^{-T}\bigg(\begin{matrix}
    -1 \\-1
    \end{matrix}\bigg) + v_2 \otimes P^{-T}\left(\begin{matrix}
    1 \\0
    \end{matrix}\right) + v_3 \otimes P^{-T}\left(\begin{matrix}
    0 \\1
    \end{matrix}\right)\\
    & = \frac{1}{\det P}\left(
    v_1 \otimes \bigg(\begin{matrix}P_{21}-P_{22}\\P_{12}-P_{11}\end{matrix}\bigg) +v_2 \otimes \bigg(\begin{matrix}P_{22}\\-P_{12}\end{matrix}\bigg)  +v_3 \otimes \bigg(\begin{matrix}-P_{21}\\P_{11}\end{matrix}\bigg)
    \right),
    \end{aligned}
\end{equation}
where we denoted $P_{ij}:=(Pe_j)\cdot e_i$. 
Clearly, if $v_i=w_i+t$ for $i=1,2,3$, with a common translation vector $t\in \R^2$, then one can use $w_i$ instead of $v_i$ in \eqref{const-grad}.
Similarly, if $v_{i} = Ap_{i} + w_{i}$ for $i = 1, 2, 3$, with $A \in \R^{2 \times 2}$ being a common matrix, then one can use $w_{i}$ instead of $v_{i}$ in \eqref{const-grad} to compute $\nabla I_{\Delta} - A$.

We are now ready to define a piecewise affine map, with constant gradients in the components of the partition of $Q_i^k$ illustrated in Figure \ref{Fig:squares}.

\smallskip

\noindent\textbf{Regions $Q_{1,k}^l$ and $Q_{2,k}^r$.} We define deformations in these regions as the boundary conditions in the neighbouring $\Sigma_{\pm \theta}$, that is 
\begin{align*}
I_{Q_{1,k}^l}(x)& \coloneq R_{-\theta} x, \quad x\in Q_{1,k}^l,\\
I_{Q_{2,k}^r}(x)& \coloneq R_\theta x, \quad x\in Q_{2,k}^r.
\end{align*}
Clearly $\nabla I_{Q_{1,k}^l}=R_{-\theta}$ and $\nabla I_{Q_{2,k}^r}=R_{\theta}$.

\smallskip

\noindent
\textbf{Regions $T^{a}_{i,k}$ and $T^{b}_{i,k}$.} We define deformations in these regions as the affine maps
\begin{align*}
I_{T^{a}_{i,k}}(x)& \coloneq R_i x+(-1)^{i+1} k \tau \varepsilon b_i, \quad x\in T^{a}_{i,k},\\
I_{T^{b}_{i,k}}(x)& \coloneq R_i x+(-1)^{i+1}(k+1) \tau \varepsilon b_i, \quad x\in T^{b}_{i,k}.
\end{align*}
Clearly $\nabla I_{T^{a}_{i,k}}=\nabla I_{T^{b}_{i,k}}=R_i$.
In particular the gradients of the deformations are independent of $k$, while the deformations are $k$-dependent.

\smallskip

\noindent
\textbf{Regions $\Delta_{i,k}^{a,n}$.} We define $I_{\Delta_{i,k}^{a,n}}$ as the interpolation in \eqref{interpol}, with points
\begin{align*}
    p_1 & =((-1)^{i+1}r_{n-1}, -r_{n-1}) +t_{i,k},\\
    p_2 & =((-1)^{i+1}r_{n},r_{n}) + t_{i,k},\\ 
    p_3 & =((-1)^{i+1}r_{n-1}, r_{n-1}) +t_{i,k},
\end{align*}
and corresponding values
$$
v_1 \coloneq R_{i} p_1 + (-1)^{i+1}(k+1) \tau\varepsilon b_i, \quad 
v_j:= R_{i} p_j + (-1)^{i+1} k \tau\varepsilon b_i, \, j=2,3.
$$
By \eqref{const-grad}, the gradient of the interpolation is 
\begin{align*}
  \nabla I_{\Delta_{i,k}^{a,n}} &= R_{i} + \frac{(-1)^{i+1}}{r_{n}r_{n-1}}\left((-1)^{i+1} \tau \varepsilon b_{i} \otimes \begin{pmatrix}
          r_{n-1} \\
          (-1)^{i}r_{n-1}
\end{pmatrix}\right) = R_{i} + \frac{1}{r_{n}}\tau \varepsilon b_{i} \otimes
\begin{pmatrix}
    1\\
    (-1)^{i}
\end{pmatrix},
\end{align*}
where we have used that $r_{n-1}=r_{n}/2$. 

\smallskip

\noindent
\textbf{Region $\Delta_{i,k}^{b,n}$.} The interpolation  $I_{\Delta_{i,k}^{b,n}}$ is defined as in \eqref{interpol}, with points 
\begin{align*}
    p_1 & = ((-1)^{i+1} r_{n-1}, -r_{n-1}) +t_{i,k},\\
    p_2 & = ((-1)^{i+1}r_{n},-r_{n}) +t_{i,k},\\
    p_3 & = ((-1)^{i+1}r_{n}, r_{n}) +t_{i,k},
\end{align*}
and corresponding values
$$
v_j:=  R_{i}p_j+(-1)^{i+1} (k+1) \tau\varepsilon b_i, \, j=1,2,\quad  
v_3:=  R_{i}p_3 + (-1)^{i+1}k \tau\varepsilon b_i.
$$
By \eqref{const-grad}, the gradient of the interpolation is
\begin{align*}
    \nabla I_{\Delta_{i,k}^{b,n}} &= R_{i} + \frac{(-1)^{i+1}}{4r_{n-1}^{2}}\left(-(-1)^{i+1} \tau \varepsilon b_{i} \otimes 
    \begin{pmatrix}
    r_{n-1} \\
    (-1)^{i+1}r_{n-1}
    \end{pmatrix}
    \right) 
     = R_{i}+\frac{1}{2r_{n}}\, \tau\varepsilon b_i\otimes
    \begin{pmatrix}
    -1 \\
    (-1)^{i}
    \end{pmatrix}.
\end{align*}

\smallskip

\noindent 
\textbf{The definition of the piecewise constant strain.}
We define the strain $\beta:\Omega \to \R^{2\times 2}$ as the piecewise constant map whose constant values are the gradients of the affine deformations of the regions of $\Omega$ defined above. In $D_{0} + t_{i,k}$, where no deformation has been defined, we simply put $\beta = \Id$.
\begin{equation}\label{strain-beta}
\beta \coloneq 
\begin{cases}
\smallskip
\Id & \text{in } \Omega \cap (D_{0}+t_{i,k}),\, k=1,\dots,N_i,\, i=1,2,\\
\smallskip
\nabla I_{\Delta^{z,n}_{i,k}}
&\text{in } \Omega \cap \Delta^{z,n}_{i,k}, \,z\in\{a,b\},\, n=1,\dots,\bar n_i, \, k=1,\dots,N_i, \, i=1,2,\\
\medskip
R_\theta &\text{in } \Sigma_\theta \cup (\Omega \cap \bigcup_{k=1}^{N_2}(Q_{2,k}^r\cup T_{2,k}^{a}\cup T_{2,k}^{b})),\\
R_{-\theta} &\text{in } \Sigma_{-\theta} \cup (\Omega \cap \bigcup_{k=1}^{N_1}(Q_{1,k}^l\cup T_{1,k}^{a}\cup T_{1,k}^{b})).
\end{cases}
\end{equation}
Note that so far the parameters $\ell_i$ (or equivalently $r_{\bar n_{i}}$), namely the width of the strips $\Sigma_i$, have not been fixed. This freedom will be needed to ensure the admissibility of the strain in Section \ref{sect:constraints}.

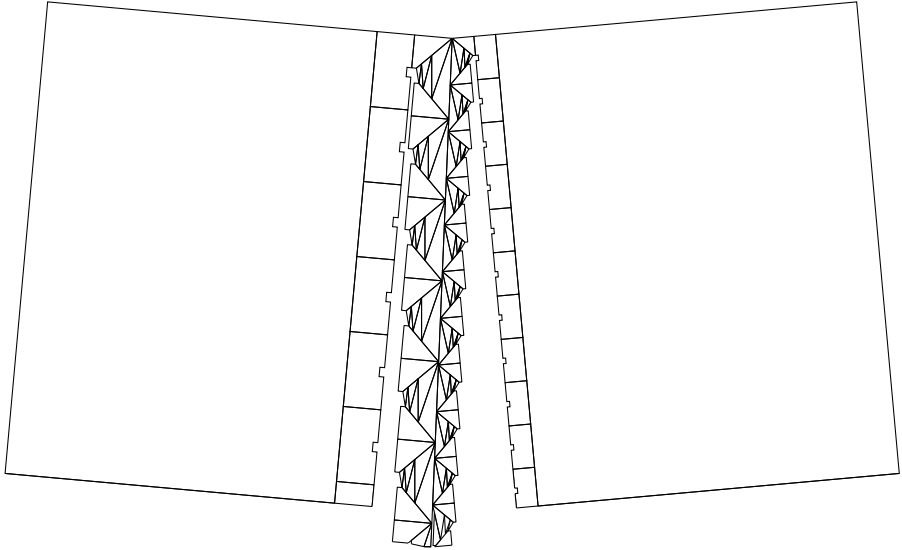
\begin{figure}[ht]
    \centering
     \begin{tikzpicture}[scale=.18]
        \tikzmath{%
            real \bx,\by,\tx,\ty,\rn,\thetaangle,\etaangle,\varphiangle,\taul,\epsl;
            \thetaangle = deg(0.09);
            \etaangle = deg(pi / 6);
            \varphiangle = deg(- pi / 3);
            \taul = 1;
            \epsl = 0.5;
            \bx1 = cos(\varphiangle);
            \by1 = sin(\varphiangle);
            \bx2 = cos(\etaangle);
            \by2 = sin(\etaangle);
            \rn1 = - (\taul * \epsl * sin(\varphiangle - \etaangle)) / (4 * sin(\thetaangle) * sin(\etaangle);
            \rn2 = (\taul * \epsl * sin(\varphiangle - \etaangle)) / (4 * sin(\thetaangle) * sin(\varphiangle));
            \tx1 = -\rn1;
            \ty1 = -\rn1;
            \tx2 = \rn2;
            \ty2 = -\rn2;
            function erre(\n,\rn){%
                real \r;
                \r = 2^(\n-3) * \rn;
                return \r;
            };
            function Rx(\x,\y,\a){%
                real \rot;
                \rot = cos(\a) * \x - sin(\a) * \y;
                return \rot;
            };
            function Ry(\x,\y,\a){%
                real \rot;
                \rot = sin(\a) * \x + cos(\a) * \y;
                return \rot;
            };
            real \px, \py;
            real \vx, \vy;
            real \pcx, \pcy;
            real \vcx, \vcy;
            \px1 = 2 * \rn2;
            \py1 = 0;
            \px2 = 30;
            \py2 = 0;
            \px3 = 30;
            \py3 = - 35;
            \px4 = 2 * \rn2;
            \py4 = - 35;
            \vx1 = Rx(\px1,\py1,\thetaangle); 
            \vy1 = Ry(\px1,\py1,\thetaangle);
            \vx2 = Rx(\px2,\py2,\thetaangle); 
            \vy2 = Ry(\px2,\py2,\thetaangle);
            \vx3 = Rx(\px3,\py3,\thetaangle); 
            \vy3 = Ry(\px3,\py3,\thetaangle);
            \vx4 = Rx(\px4,\py4,\thetaangle); 
            \vy4 = Ry(\px4,\py4,\thetaangle);
            {\draw (\vx1,\vy1) -- (\vx2,\vy2) -- (\vx3,\vy3) -- (\vx4,\vy4) -- cycle;};
            \px1 = -2 * \rn1;
            \py1 = 0;
            \px2 = - 30;
            \py2 = 0;
            \px3 = - 30;
            \py3 = - 35;
            \px4 = -2 * \rn1;
            \py4 = - 35;
            \vx1 = Rx(\px1,\py1,-\thetaangle); 
            \vy1 = Ry(\px1,\py1,-\thetaangle);
            \vx2 = Rx(\px2,\py2,-\thetaangle); 
            \vy2 = Ry(\px2,\py2,-\thetaangle);
            \vx3 = Rx(\px3,\py3,-\thetaangle); 
            \vy3 = Ry(\px3,\py3,-\thetaangle);
            \vx4 = Rx(\px4,\py4,-\thetaangle); 
            \vy4 = Ry(\px4,\py4,-\thetaangle);
            {\draw (\vx1,\vy1) -- (\vx2,\vy2) -- (\vx3,\vy3) -- (\vx4,\vy4) -- cycle;};
            int \k;
            \pcx1 = \tx1;
            \pcy1 = 0;
            \pcx2 = - 30;
            \pcy2 = 0;
            \pcx3 = - 30;
            \pcy3 = - 35;
            \pcx4 = \tx1;
            \pcy4 = - 35;
            \vcx1 = Rx(\pcx1,\pcy1,-\thetaangle); 
            \vcy1 = Ry(\pcx1,\pcy1,-\thetaangle);
            \vcx2 = Rx(\pcx2,\pcy2,-\thetaangle); 
            \vcy2 = Ry(\pcx2,\pcy2,-\thetaangle);
            \vcx3 = Rx(\pcx3,\pcy3,-\thetaangle); 
            \vcy3 = Ry(\pcx3,\pcy3,-\thetaangle);
            \vcx4 = Rx(\pcx4,\pcy4,-\thetaangle); 
            \vcy4 = Ry(\pcx4,\pcy4,-\thetaangle);
            for \k in {0,...,6} {
                \px1 = \tx1;
                \py1 = (2 * \k + 1) * \ty1 + \rn1;
                \px2 = \tx1 - \rn1;
                \py2 = (2 * \k + 1) * \ty1 + \rn1;
                \px3 = \tx1 - \rn1;
                \py3 = (2 * \k + 1) * \ty1 - \rn1;
                \px4 = \tx1;
                \py4 = (2 * \k + 1) * \ty1 - \rn1;
                \px5 = \tx1;
                \py5 = (2 * \k + 1) * \ty1 - erre(0,\rn1);
                \px6 = \tx1 - erre(0, \rn1);
                \py6 = (2 * \k + 1) * \ty1 - erre(0,\rn1);
                \px7 = \tx1 - erre(0, \rn1);
                \py7 = (2 * \k + 1) * \ty1 + erre(0,\rn1);
                \px8 = \tx1;
                \py8 = (2 * \k + 1) * \ty1 + erre(0,\rn1);
                \vx1 = Rx(\px1,\py1,-\thetaangle); 
                \vy1 = Ry(\px1,\py1,-\thetaangle);
                \vx2 = Rx(\px2,\py2,-\thetaangle); 
                \vy2 = Ry(\px2,\py2,-\thetaangle);
                \vx3 = Rx(\px3,\py3,-\thetaangle); 
                \vy3 = Ry(\px3,\py3,-\thetaangle);
                \vx4 = Rx(\px4,\py4,-\thetaangle); 
                \vy4 = Ry(\px4,\py4,-\thetaangle);
                \vx5 = Rx(\px5,\py5,-\thetaangle); 
                \vy5 = Ry(\px5,\py5,-\thetaangle);
                \vx6 = Rx(\px6,\py6,-\thetaangle); 
                \vy6 = Ry(\px6,\py6,-\thetaangle);
                \vx7 = Rx(\px7,\py7,-\thetaangle); 
                \vy7 = Ry(\px7,\py7,-\thetaangle);
                \vx8 = Rx(\px8,\py8,-\thetaangle); 
                \vy8 = Ry(\px8,\py8,-\thetaangle);
                if \k == 6 then {
                    {%
                        \begin{scope}
                            \clip (\vcx1,\vcy1) -- (\vcx2,\vcy2) -- (\vcx3,\vcy3) -- (\vcx4,\vcy4) -- cycle;
                            \draw (\vx1,\vy1) -- (\vx2,\vy2) -- (\vx3,\vy3) -- (\vx4,\vy4) -- (\vx5,\vy5) -- (\vx6,\vy6) -- (\vx7,\vy7) -- (\vx8,\vy8) -- cycle;
                        \end{scope}
                        \draw (\vcx3,\vcy3) -- (\vcx4,\vcy4);
                        \draw (\vcx4,\vcy4) -- (\vx1,\vy1);
                    };
                } else {
                    {
                        \draw (\vx1,\vy1) -- (\vx2,\vy2) -- (\vx3,\vy3) -- (\vx4,\vy4) -- (\vx5,\vy5) -- (\vx6,\vy6) -- (\vx7,\vy7) -- (\vx8,\vy8) -- cycle;
                    };
                };
            };
            \pcx1 = \tx2;
            \pcy1 = 0;
            \pcx2 = 30;
            \pcy2 = 0;
            \pcx3 = 30;
            \pcy3 = - 35;
            \pcx4 = \tx2;
            \pcy4 = - 35;
            \vcx1 = Rx(\pcx1,\pcy1,\thetaangle); 
            \vcy1 = Ry(\pcx1,\pcy1,\thetaangle);
            \vcx2 = Rx(\pcx2,\pcy2,\thetaangle); 
            \vcy2 = Ry(\pcx2,\pcy2,\thetaangle);
            \vcx3 = Rx(\pcx3,\pcy3,\thetaangle); 
            \vcy3 = Ry(\pcx3,\pcy3,\thetaangle);
            \vcx4 = Rx(\pcx4,\pcy4,\thetaangle); 
            \vcy4 = Ry(\pcx4,\pcy4,\thetaangle);
            for \k in {0,...,10} {
                \px1 = \tx2;
                \py1 = (2 * \k + 1) * \ty2 + \rn2;
                \px2 = \tx2 + \rn2;
                \py2 = (2 * \k + 1) * \ty2 + \rn2;
                \px3 = \tx2 + \rn2;
                \py3 = (2 * \k + 1) * \ty2 - \rn2;
                \px4 = \tx2;
                \py4 = (2 * \k + 1) * \ty2 - \rn2;
                \px5 = \tx2;
                \py5 = (2 * \k + 1) * \ty2 - erre(0,\rn2);
                \px6 = \tx2 + erre(0, \rn2);
                \py6 = (2 * \k + 1) * \ty2 - erre(0,\rn2);
                \px7 = \tx2 + erre(0, \rn2);
                \py7 = (2 * \k + 1) * \ty2 + erre(0,\rn2);
                \px8 = \tx2;
                \py8 = (2 * \k + 1) * \ty2 + erre(0,\rn2);
                \vx1 = Rx(\px1,\py1,\thetaangle); 
                \vy1 = Ry(\px1,\py1,\thetaangle);
                \vx2 = Rx(\px2,\py2,\thetaangle); 
                \vy2 = Ry(\px2,\py2,\thetaangle);
                \vx3 = Rx(\px3,\py3,\thetaangle); 
                \vy3 = Ry(\px3,\py3,\thetaangle);
                \vx4 = Rx(\px4,\py4,\thetaangle); 
                \vy4 = Ry(\px4,\py4,\thetaangle);
                \vx5 = Rx(\px5,\py5,\thetaangle); 
                \vy5 = Ry(\px5,\py5,\thetaangle);
                \vx6 = Rx(\px6,\py6,\thetaangle); 
                \vy6 = Ry(\px6,\py6,\thetaangle);
                \vx7 = Rx(\px7,\py7,\thetaangle); 
                \vy7 = Ry(\px7,\py7,\thetaangle);
                \vx8 = Rx(\px8,\py8,\thetaangle); 
                \vy8 = Ry(\px8,\py8,\thetaangle);
                if \k == 10 then {
                    {%
                        \begin{scope}
                            \clip (\vcx1,\vcy1) -- (\vcx2,\vcy2) -- (\vcx3,\vcy3) -- (\vcx4,\vcy4) -- cycle;
                            \draw (\vx1,\vy1) -- (\vx2,\vy2) -- (\vx3,\vy3) -- (\vx4,\vy4) -- (\vx5,\vy5) -- (\vx6,\vy6) -- (\vx7,\vy7) -- (\vx8,\vy8) --cycle;
                        \end{scope}
                        \draw (\vcx3,\vcy3) -- (\vcx4,\vcy4);
                        \draw (\vcx4,\vcy4) -- (\vx5,\vy5);
                        \draw (\vx8,\vy8) -- (\vx1,\vy1);
                    };
                } else {
                    {
                        \draw (\vx1,\vy1) -- (\vx2,\vy2) -- (\vx3,\vy3) -- (\vx4,\vy4) -- (\vx5,\vy5) -- (\vx6,\vy6) -- (\vx7,\vy7) -- (\vx8,\vy8) --cycle;
                    };
                };
            };
            \pcx1 = 0;
            \pcy1 = 0;
            \pcx2 = -30;
            \pcy2 = 0;
            \pcx3 = -30;
            \pcy3 = - 35;
            \pcx4 = 0;
            \pcy4 = - 35;
            \vcx1 = Rx(\pcx1,\pcy1,-\thetaangle) + (6) * \taul * \epsl * \bx1; 
            \vcy1 = Ry(\pcx1,\pcy1,-\thetaangle) + (6) * \taul * \epsl * \by1;
            \vcx2 = Rx(\pcx2,\pcy2,-\thetaangle) + (6) * \taul * \epsl * \bx1; 
            \vcy2 = Ry(\pcx2,\pcy2,-\thetaangle) + (6) * \taul * \epsl * \by1;
            \vcx3 = Rx(\pcx3,\pcy3,-\thetaangle) + (6) * \taul * \epsl * \bx1; 
            \vcy3 = Ry(\pcx3,\pcy3,-\thetaangle) + (6) * \taul * \epsl * \by1;
            \vcx4 = Rx(\pcx4,\pcy4,-\thetaangle) + (6) * \taul * \epsl * \bx1; 
            \vcy4 = Ry(\pcx4,\pcy4,-\thetaangle) + (6) * \taul * \epsl * \by1;
            for \k in {0,...,6} {
                \px1 = \tx1;
                \py1 = (2 * \k + 1) * \ty1 + \rn1;
                \px2 = \tx1 + \rn1;
                \py2 = (2 * \k + 1) * \ty1 + \rn1;
                \px3 = \tx1 + erre(0,\rn1);
                \py3 = (2 * \k + 1) * \ty1 + erre(0,\rn1);
                \px4 = \tx1;
                \py4 = (2 * \k + 1) * \ty1 + erre(0,\rn1);
                \vx1 = Rx(\px1,\py1,-\thetaangle) + \k * \taul * \epsl * \bx1; 
                \vy1 = Ry(\px1,\py1,-\thetaangle) + \k * \taul * \epsl * \by1;
                \vx2 = Rx(\px2,\py2,-\thetaangle) + \k * \taul * \epsl * \bx1; 
                \vy2 = Ry(\px2,\py2,-\thetaangle) + \k * \taul * \epsl * \by1;
                \vx3 = Rx(\px3,\py3,-\thetaangle) + \k * \taul * \epsl * \bx1; 
                \vy3 = Ry(\px3,\py3,-\thetaangle) + \k * \taul * \epsl * \by1;
                \vx4 = Rx(\px4,\py4,-\thetaangle) + \k * \taul * \epsl * \bx1; 
                \vy4 = Ry(\px4,\py4,-\thetaangle) + \k * \taul * \epsl * \by1;
                if \k == 6 then {
                    {%
                        \begin{scope}
                            \draw[name path global= clipper, transparent] (\vcx1,\vcy1) -- (\vcx2,\vcy2) -- (\vcx3,\vcy3) -- (\vcx4,\vcy4) -- cycle;
                            \clip (\vcx1,\vcy1) -- (\vcx2,\vcy2) -- (\vcx3,\vcy3) -- (\vcx4,\vcy4) -- cycle;
                            \draw[name path global = last, line join = bevel] (\vx1,\vy1) -- (\vx2,\vy2) -- (\vx3,\vy3) -- (\vx4,\vy4) -- cycle;      
                        \end{scope}
                        \draw[intersection segments={of=clipper and last,sequence = L2}];
                    };
                } else {
                    {\draw[line join = bevel] (\vx1,\vy1) -- (\vx2,\vy2) -- (\vx3,\vy3) -- (\vx4,\vy4) -- cycle;};
                };
            };
            for \k in {0,...,5} {
                \px1 = \tx1;
                \py1 = (2 * \k + 1) * \ty1 - \rn1;
                \px2 = \tx1 + \rn1;
                \py2 = (2 * \k + 1) * \ty1 - \rn1;
                \px3 = \tx1 + erre(0,\rn1);
                \py3 = (2 * \k + 1) * \ty1 - erre(0,\rn1);
                \px4 = \tx1;
                \py4 = (2 * \k + 1) * \ty1 - erre(0,\rn1);
                \vx1 = Rx(\px1,\py1,-\thetaangle) + (\k + 1) * \taul * \epsl * \bx1; 
                \vy1 = Ry(\px1,\py1,-\thetaangle) + (\k + 1) * \taul * \epsl * \by1;
                \vx2 = Rx(\px2,\py2,-\thetaangle) + (\k + 1) * \taul * \epsl * \bx1; 
                \vy2 = Ry(\px2,\py2,-\thetaangle) + (\k + 1) * \taul * \epsl * \by1;
                \vx3 = Rx(\px3,\py3,-\thetaangle) + (\k + 1) * \taul * \epsl * \bx1; 
                \vy3 = Ry(\px3,\py3,-\thetaangle) + (\k + 1) * \taul * \epsl * \by1;
                \vx4 = Rx(\px4,\py4,-\thetaangle) + (\k + 1) * \taul * \epsl * \bx1; 
                \vy4 = Ry(\px4,\py4,-\thetaangle) + (\k + 1) * \taul * \epsl * \by1;
                {\draw[line join = bevel] (\vx1,\vy1) -- (\vx2,\vy2) -- (\vx3,\vy3) -- (\vx4,\vy4) -- cycle;};
            };
            for \k in {0,...,10} {
                \px1 = \tx2;
                \py1 = (2 * \k + 1) * \ty2 + \rn2;
                \px2 = \tx2 - \rn2;
                \py2 = (2 * \k + 1) * \ty2 + \rn2;
                \px3 = \tx2 - erre(0,\rn2);
                \py3 = (2 * \k + 1) * \ty2 + erre(0,\rn2);
                \px4 = \tx2;
                \py4 = (2 * \k + 1) * \ty2 + erre(0,\rn2);
                \vx1 = Rx(\px1,\py1,\thetaangle) - \k * \taul * \epsl * \bx2; 
                \vy1 = Ry(\px1,\py1,\thetaangle) - \k * \taul * \epsl * \by2;
                \vx2 = Rx(\px2,\py2,\thetaangle) - \k * \taul * \epsl * \bx2; 
                \vy2 = Ry(\px2,\py2,\thetaangle) - \k * \taul * \epsl * \by2;
                \vx3 = Rx(\px3,\py3,\thetaangle) - \k * \taul * \epsl * \bx2; 
                \vy3 = Ry(\px3,\py3,\thetaangle) - \k * \taul * \epsl * \by2;
                \vx4 = Rx(\px4,\py4,\thetaangle) - \k * \taul * \epsl * \bx2; 
                \vy4 = Ry(\px4,\py4,\thetaangle) - \k * \taul * \epsl * \by2;
                {\draw[line join = bevel] (\vx1,\vy1) -- (\vx2,\vy2) -- (\vx3,\vy3) -- (\vx4,\vy4) -- cycle;};
            };
            \pcx1 = 0;
            \pcy1 = 0;
            \pcx2 = 30;
            \pcy2 = 0;
            \pcx3 = 30;
            \pcy3 = - 35;
            \pcx4 = 0;
            \pcy4 = - 35;
            \vcx1 = Rx(\pcx1,\pcy1,\thetaangle) - (11) * \taul * \epsl * \bx2; 
            \vcy1 = Ry(\pcx1,\pcy1,\thetaangle) - (11) * \taul * \epsl * \by2;
            \vcx2 = Rx(\pcx2,\pcy2,\thetaangle) - (11) * \taul * \epsl * \bx2; 
            \vcy2 = Ry(\pcx2,\pcy2,\thetaangle) - (11) * \taul * \epsl * \by2;
            \vcx3 = Rx(\pcx3,\pcy3,\thetaangle) - (11) * \taul * \epsl * \bx2; 
            \vcy3 = Ry(\pcx3,\pcy3,\thetaangle) - (11) * \taul * \epsl * \by2;
            \vcx4 = Rx(\pcx4,\pcy4,\thetaangle) - (11) * \taul * \epsl * \bx2; 
            \vcy4 = Ry(\pcx4,\pcy4,\thetaangle) - (11) * \taul * \epsl * \by2;
            for \k in {0,...,10} {
                \px1 = \tx2;
                \py1 = (2 * \k + 1) * \ty2 - \rn2;
                \px2 = \tx2 - \rn2;
                \py2 = (2 * \k + 1) * \ty2 - \rn2;
                \px3 = \tx2 - erre(0,\rn2);
                \py3 = (2 * \k + 1) * \ty2 - erre(0,\rn2);
                \px4 = \tx2; \py4 = (2 * \k + 1) * \ty2 - erre(0,\rn2);
                \vx1 = Rx(\px1,\py1,\thetaangle) - (\k + 1) * \taul * \epsl * \bx2; 
                \vy1 = Ry(\px1,\py1,\thetaangle) - (\k + 1) * \taul * \epsl * \by2;
                \vx2 = Rx(\px2,\py2,\thetaangle) - (\k + 1) * \taul * \epsl * \bx2; 
                \vy2 = Ry(\px2,\py2,\thetaangle) - (\k + 1) * \taul * \epsl * \by2;
                \vx3 = Rx(\px3,\py3,\thetaangle) - (\k + 1) * \taul * \epsl * \bx2; 
                \vy3 = Ry(\px3,\py3,\thetaangle) - (\k + 1) * \taul * \epsl * \by2;
                \vx4 = Rx(\px4,\py4,\thetaangle) - (\k + 1) * \taul * \epsl * \bx2; 
                \vy4 = Ry(\px4,\py4,\thetaangle) - (\k + 1) * \taul * \epsl * \by2;
                if \k == 10 then {
                    {%
                        \begin{scope}
                            \draw[name path global= clipper, transparent] (\vcx1,\vcy1) -- (\vcx2,\vcy2) -- (\vcx3,\vcy3) -- (\vcx4,\vcy4) -- cycle;
                            \clip (\vcx1,\vcy1) -- (\vcx2,\vcy2) -- (\vcx3,\vcy3) -- (\vcx4,\vcy4) -- cycle;
                            \draw[name path global= last, line join = bevel] (\vx1,\vy1) -- (\vx2,\vy2) -- (\vx3,\vy3) -- (\vx4,\vy4) -- cycle;
                        \end{scope}
                        \draw[intersection segments={of=clipper and last,sequence = L2}];
                    };
                } else {
                    {\draw[line join = bevel] (\vx1,\vy1) -- (\vx2,\vy2) -- (\vx3,\vy3) -- (\vx4,\vy4) -- cycle;};
                };
            };
            int \n;
            for \n in {1,2,3}{
                \pcx1 = 0;
                \pcy1 = 0;
                \pcx2 = -30;
                \pcy2 = 0;
                \pcx3 = -30;
                \pcy3 = - 35;
                \pcx4 = 0;
                \pcy4 = - 35;
                \vcx1 = Rx(\pcx1,\pcy1,-\thetaangle) + (\taul * \epsl / erre(\n, \rn1)) * \bx1 * (\pcx1 - \pcy1) - (\taul * \epsl / erre(\n, \rn1)) * \bx1 * (-\rn1 + 13 * \rn1) + 6 * \taul * \epsl * \bx1; 
                \vcy1 = Ry(\pcx1,\pcy1,-\thetaangle) + (\taul * \epsl / erre(\n, \rn1)) * \by1 * (\pcx1 - \pcy1) - (\taul * \epsl / erre(\n, \rn1)) * \by1 * (-\rn1 + 13 * \rn1) + 6 * \taul * \epsl * \by1;
                \vcx2 = Rx(\pcx2,\pcy2,-\thetaangle) + (\taul * \epsl / erre(\n, \rn1)) * \bx1 * (\pcx2 - \pcy2) - (\taul * \epsl / erre(\n, \rn1)) * \bx1 * (-\rn1 + 13 * \rn1) + 6 * \taul * \epsl * \bx1;
                \vcy2 = Ry(\pcx2,\pcy2,-\thetaangle) + (\taul * \epsl / erre(\n, \rn1)) * \by1 * (\pcx2 - \pcy2) - (\taul * \epsl / erre(\n, \rn1)) * \by1 * (-\rn1 + 13 * \rn1) + 6 * \taul * \epsl * \by1;
                \vcx3 = Rx(\pcx3,\pcy3,-\thetaangle) + (\taul * \epsl / erre(\n, \rn1)) * \bx1 * (\pcx3 - \pcy3) - (\taul * \epsl / erre(\n, \rn1)) * \bx1 * (-\rn1 + 13 * \rn1) + 6 * \taul * \epsl * \bx1;
                \vcy3 = Ry(\pcx3,\pcy3,-\thetaangle) + (\taul * \epsl / erre(\n, \rn1)) * \by1 * (\pcx3 - \pcy3) - (\taul * \epsl / erre(\n, \rn1)) * \by1 * (-\rn1 + 13 * \rn1) + 6 * \taul * \epsl * \by1;
                \vcx4 = Rx(\pcx4,\pcy4,-\thetaangle) + (\taul * \epsl / erre(\n, \rn1)) * \bx1 * (\pcx4 - \pcy4) - (\taul * \epsl / erre(\n, \rn1)) * \bx1 * (-\rn1 + 13 * \rn1) + 6 * \taul * \epsl * \bx1;
                \vcy4 = Ry(\pcx4,\pcy4,-\thetaangle) + (\taul * \epsl / erre(\n, \rn1)) * \by1 * (\pcx4 - \pcy4) - (\taul * \epsl / erre(\n, \rn1)) * \by1 * (-\rn1 + 13 * \rn1) + 6 * \taul * \epsl * \by1;
                for \k in {0,...,6} {
                    \px1 = \tx1 + erre(\n - 1, \rn1);
                    \py1 = (2 * \k + 1) * \ty1 + erre(\n - 1, \rn1);
                    \px2 = \tx1 + erre(\n, \rn1);
                    \py2 = (2 * \k + 1) * \ty1 + erre(\n, \rn1);
                    \px3 = \tx1 + erre(\n-1, \rn1);
                    \py3 = (2 * \k + 1) * \ty1 - erre(\n-1, \rn1);
                    \vx1 = Rx(\px1,\py1,-\thetaangle) + \k * \taul * \epsl * \bx1; 
                    \vy1 = Ry(\px1,\py1,-\thetaangle) + \k * \taul * \epsl * \by1; 
                    \vx2 = Rx(\px2,\py2,-\thetaangle) + \k * \taul * \epsl * \bx1; 
                    \vy2 = Ry(\px2,\py2,-\thetaangle) + \k * \taul * \epsl * \by1; 
                    \vx3 = Rx(\px3,\py3,-\thetaangle) + (\k + 1) * \taul * \epsl * \bx1; 
                    \vy3 = Ry(\px3,\py3,-\thetaangle) + (\k + 1) * \taul * \epsl * \by1; 
                    if \k == 6 then {
                        if \n == 3 then {
                            {%
                                \begin{scope}
                                    \draw[name path global = clipper, transparent] (\vcx1,\vcy1) -- (\vcx2,\vcy2) -- (\vcx3,\vcy3) -- (\vcx4,\vcy4) -- cycle;
                                    \clip (\vcx1,\vcy1) -- (\vcx2,\vcy2) -- (\vcx3,\vcy3) -- (\vcx4,\vcy4) -- cycle;
                                    \draw[name path global = last, line join = bevel] (\vx1,\vy1) -- (\vx2,\vy2) -- (\vx3,\vy3) -- cycle;
                                \end{scope}
                                \draw[intersection segments={of=clipper and last,sequence = L2}];
                            };
                        };
                    } else {
                        {\draw[line join = bevel] (\vx1,\vy1) -- (\vx2,\vy2) -- (\vx3,\vy3) -- cycle;};
                    };
                };
            };
            int \n;
            for \n in {1,2,3}{
                \pcx1 = 0;
                \pcy1 = 0;
                \pcx2 = -30;
                \pcy2 = 0;
                \pcx3 = -30;
                \pcy3 = - 35;
                \pcx4 = 0;
                \pcy4 = - 35;
                \vcx1 = Rx(\pcx1,\pcy1,-\thetaangle) - (\taul * \epsl / (2 * erre(\n, \rn1))) * \bx1 * (\pcx1 + \pcy1) + (\taul * \epsl / (2 * erre(\n, \rn1))) * \bx1 * (-\rn1 - 13 * \rn1) + 7 * \taul * \epsl * \bx1; 
                \vcy1 = Ry(\pcx1,\pcy1,-\thetaangle) - (\taul * \epsl / (2 * erre(\n, \rn1))) * \by1 * (\pcx1 + \pcy1) + (\taul * \epsl / (2 * erre(\n, \rn1))) * \by1 * (-\rn1 - 13 * \rn1) + 7 * \taul * \epsl * \by1;
                \vcx2 = Rx(\pcx2,\pcy2,-\thetaangle) - (\taul * \epsl / (2 * erre(\n, \rn1))) * \bx1 * (\pcx2 + \pcy2) + (\taul * \epsl / (2 * erre(\n, \rn1))) * \bx1 * (-\rn1 - 13 * \rn1) + 7 * \taul * \epsl * \bx1;
                \vcy2 = Ry(\pcx2,\pcy2,-\thetaangle) - (\taul * \epsl / (2 * erre(\n, \rn1))) * \by1 * (\pcx2 + \pcy2) + (\taul * \epsl / (2 * erre(\n, \rn1))) * \by1 * (-\rn1 - 13 * \rn1) + 7 * \taul * \epsl * \by1;
                \vcx3 = Rx(\pcx3,\pcy3,-\thetaangle) - (\taul * \epsl / (2 * erre(\n, \rn1))) * \bx1 * (\pcx3 + \pcy3) + (\taul * \epsl / (2 * erre(\n, \rn1))) * \bx1 * (-\rn1 - 13 * \rn1) + 7 * \taul * \epsl * \bx1;
                \vcy3 = Ry(\pcx3,\pcy3,-\thetaangle) - (\taul * \epsl / (2 * erre(\n, \rn1))) * \by1 * (\pcx3 + \pcy3) + (\taul * \epsl / (2 * erre(\n, \rn1))) * \by1 * (-\rn1 - 13 * \rn1) + 7 * \taul * \epsl * \by1;
                \vcx4 = Rx(\pcx4,\pcy4,-\thetaangle) - (\taul * \epsl / (2 * erre(\n, \rn1))) * \bx1 * (\pcx4 + \pcy4) + (\taul * \epsl / (2 * erre(\n, \rn1))) * \bx1 * (-\rn1 - 13 * \rn1) + 7 * \taul * \epsl * \bx1;
                \vcy4 = Ry(\pcx4,\pcy4,-\thetaangle) - (\taul * \epsl / (2 * erre(\n, \rn1))) * \by1 * (\pcx4 + \pcy4) + (\taul * \epsl / (2 * erre(\n, \rn1))) * \by1 * (-\rn1 - 13 * \rn1) + 7 * \taul * \epsl * \by1;
                for \k in {0,...,6} {
                    \px1 = \tx1 + erre(\n, \rn1);
                    \py1 = (2 * \k + 1) * \ty1 + erre(\n, \rn1);
                    \px2 = \tx1 + erre(\n, \rn1);
                    \py2 = (2 * \k + 1) * \ty1 - erre(\n, \rn1);
                    \px3 = \tx1 + erre(\n-1, \rn1);
                    \py3 = (2 * \k + 1) * \ty1 - erre(\n-1, \rn1);
                    \vx1 = Rx(\px1,\py1,-\thetaangle) - (\taul * \epsl / (2 * erre(\n, \rn1))) * (\bx1 * erre(\n, \rn1) + \bx1 * erre(\n, \rn1)) + (\k + 1) * \taul * \epsl * \bx1; 
                    \vy1 = Ry(\px1,\py1,-\thetaangle) - (\taul * \epsl / (2 * erre(\n, \rn1))) * (\by1 * erre(\n, \rn1) + \by1 * erre(\n, \rn1)) + (\k + 1) * \taul * \epsl * \by1; 
                    \vx2 = Rx(\px2,\py2,-\thetaangle) - (\taul * \epsl / (2 * erre(\n, \rn1))) * (\bx1 * erre(\n, \rn1) - \bx1 * erre(\n, \rn1)) + (\k + 1) * \taul * \epsl * \bx1; 
                    \vy2 = Ry(\px2,\py2,-\thetaangle) - (\taul * \epsl / (2 * erre(\n, \rn1))) * (\by1 * erre(\n,\rn1) - \by1 * erre(\n, \rn1)) + (\k + 1) * \taul * \epsl * \by1; 
                    \vx3 = Rx(\px3,\py3,-\thetaangle) - (\taul * \epsl / (2 * erre(\n, \rn1))) * (\bx1 * erre(\n-1,\rn1) - \bx1 * erre(\n-1,\rn1)) + (\k + 1) * \taul * \epsl * \bx1; 
                    \vy3 = Ry(\px3,\py3,-\thetaangle) - (\taul * \epsl / (2 * erre(\n, \rn1))) * (\by1 * erre(\n-1,\rn1) - \by1 * erre(\n-1,\rn1)) + (\k + 1) * \taul * \epsl * \by1; 
                    if \k == 6 then {
                        if \n == 3 then {
                            {%
                                \begin{scope}
                                    \draw[name path global = clipper, transparent] (\vcx1,\vcy1) -- (\vcx2,\vcy2) -- (\vcx3,\vcy3) -- (\vcx4,\vcy4) -- cycle;
                                    \clip (\vcx1,\vcy1) -- (\vcx2,\vcy2) -- (\vcx3,\vcy3) -- (\vcx4,\vcy4) -- cycle;
                                    \draw[name path global = last, line join = bevel] (\vx1,\vy1) -- (\vx2,\vy2) -- (\vx3,\vy3) -- cycle;
                                \end{scope}
                                \draw[intersection segments={of=clipper and last,sequence = {L2 -- L3}}];
                            };
                        };
                    } else {
                        {\draw[line join = bevel] (\vx1,\vy1) -- (\vx2,\vy2) -- (\vx3,\vy3) -- cycle;}; };
                    };
            };
            int \n;
            for \k in {0,...,10} {
                for \n in {1,2,3}{
                    \px1 = \tx2 - erre(\n - 1, \rn2);
                    \py1 = (2 * \k + 1) * \ty2 + erre(\n - 1, \rn2);
                    \px2 = \tx2 - erre(\n, \rn2);
                    \py2 = (2 * \k + 1) * \ty2 + erre(\n, \rn2);
                    \px3 = \tx2 - erre(\n-1, \rn2);
                    \py3 = (2 * \k + 1) * \ty2 - erre(\n-1, \rn2);
                    \vx1 = Rx(\px1,\py1,\thetaangle) + (\taul * \epsl / erre(\n, \rn2)) * (- \bx2 * erre(\n-1, \rn2) + \bx2 * erre(\n - 1, \rn2)) - \k * \taul * \epsl * \bx2; 
                    \vy1 = Ry(\px1,\py1,\thetaangle) + (\taul * \epsl / erre(\n, \rn2)) * (- \by2 * erre(\n -1, \rn2) + \by2 * erre(\n - 1, \rn2)) - \k * \taul * \epsl * \by2; 
                    \vx2 = Rx(\px2,\py2,\thetaangle) + (\taul * \epsl / erre(\n, \rn2)) * (- \bx2 * erre(\n, \rn2) + \bx2 * erre(\n, \rn2)) - \k * \taul * \epsl * \bx2; 
                    \vy2 = Ry(\px2,\py2,\thetaangle) + (\taul * \epsl / erre(\n, \rn2)) * (- \by2 * erre(\n,\rn2) + \by2 * erre(\n, \rn2)) - \k * \taul * \epsl * \by2; 
                    \vx3 = Rx(\px3,\py3,\thetaangle) + (\taul * \epsl / erre(\n, \rn2)) * (- \bx2 * erre(\n-1,\rn2) - \bx2 * erre(\n-1,\rn2)) - \k * \taul * \epsl * \bx2; 
                    \vy3 = Ry(\px3,\py3,\thetaangle) + (\taul * \epsl / erre(\n, \rn2)) * (- \by2 * erre(\n-1,\rn2) - \by2 * erre(\n-1,\rn2)) - \k * \taul * \epsl * \by2; 
                    {\draw[line join = bevel] (\vx1,\vy1) -- (\vx2,\vy2) -- (\vx3,\vy3) -- cycle;};
                };
            };
            int \n;
            for \n in {1,2,3}{
                \pcx1 = 0;
                \pcy1 = 0;
                \pcx2 = 30;
                \pcy2 = 0;
                \pcx3 = 30;
                \pcy3 = - 35;
                \pcx4 = 0;
                \pcy4 = - 35;
                \vcx1 = Rx(\pcx1,\pcy1,\thetaangle) + (\taul * \epsl / (2 * erre(\n, \rn2))) * \bx2 * (-\pcx1 + \pcy1) - (\taul * \epsl / (2 * erre(\n, \rn2))) * \bx2 * (-\rn1 - 21 * \rn2) - 11 * \taul * \epsl * \bx2; 
                \vcy1 = Ry(\pcx1,\pcy1,\thetaangle) + (\taul * \epsl / (2 * erre(\n, \rn2))) * \by2 * (-\pcx1 + \pcy1) - (\taul * \epsl / (2 * erre(\n, \rn2))) * \by2 * (-\rn1 - 21 * \rn2) - 11 * \taul * \epsl * \by2;
                \vcx2 = Rx(\pcx2,\pcy2,\thetaangle) + (\taul * \epsl / (2 * erre(\n, \rn2))) * \bx2 * (-\pcx2 + \pcy2) - (\taul * \epsl / (2 * erre(\n, \rn2))) * \bx2 * (-\rn1 - 21 * \rn2) - 11 * \taul * \epsl * \bx2;
                \vcy2 = Ry(\pcx2,\pcy2,\thetaangle) + (\taul * \epsl / (2 * erre(\n, \rn2))) * \by2 * (-\pcx2 + \pcy2) - (\taul * \epsl / (2 * erre(\n, \rn2))) * \by2 * (-\rn1 - 21 * \rn2) - 11 * \taul * \epsl * \by2;
                \vcx3 = Rx(\pcx3,\pcy3,\thetaangle) + (\taul * \epsl / (2 * erre(\n, \rn2))) * \bx2 * (-\pcx3 + \pcy3) - (\taul * \epsl / (2 * erre(\n, \rn2))) * \bx2 * (-\rn1 - 21 * \rn2) - 11 * \taul * \epsl * \bx2;
                \vcy3 = Ry(\pcx3,\pcy3,\thetaangle) + (\taul * \epsl / (2 * erre(\n, \rn2))) * \by2 * (-\pcx3 + \pcy3) - (\taul * \epsl / (2 * erre(\n, \rn2))) * \by2 * (-\rn1 - 21 * \rn2) - 11 * \taul * \epsl * \by2;
                \vcx4 = Rx(\pcx4,\pcy4,\thetaangle) + (\taul * \epsl / (2 * erre(\n, \rn2))) * \bx2 * (-\pcx4 + \pcy4) - (\taul * \epsl / (2 * erre(\n, \rn2))) * \bx2 * (-\rn1 - 21 * \rn2) - 11 * \taul * \epsl * \bx2;
                \vcy4 = Ry(\pcx4,\pcy4,\thetaangle) + (\taul * \epsl / (2 * erre(\n, \rn2))) * \by2 * (-\pcx4 + \pcy4) - (\taul * \epsl / (2 * erre(\n, \rn2))) * \by2 * (-\rn1 - 21 * \rn2) - 11 * \taul * \epsl * \by2;
                for \k in {0,...,10} {
                    \px1 = \tx2 - erre(\n, \rn2);
                    \py1 = (2 * \k + 1) * \ty2 + erre(\n, \rn2);
                    \px2 = \tx2 - erre(\n, \rn2);
                    \py2 = (2 * \k + 1) * \ty2 - erre(\n, \rn2);
                    \px3 = \tx2 - erre(\n-1, \rn2);
                    \py3 = (2 * \k + 1) * \ty2 - erre(\n-1, \rn2);
                    \vx1 = Rx(\px1,\py1,\thetaangle) - (\taul * \epsl / (2 * erre(\n, \rn2))) * (- \bx2 * erre(\n, \rn2) - \bx2 * erre(\n, \rn2)) - (\k + 1) * \taul * \epsl * \bx2; 
                    \vy1 = Ry(\px1,\py1,\thetaangle) - (\taul * \epsl / (2 * erre(\n, \rn2))) * (- \by2 * erre(\n, \rn2) - \by2 * erre(\n, \rn2)) - (\k + 1) * \taul * \epsl * \by2; 
                    \vx2 = Rx(\px2,\py2,\thetaangle) - (\taul * \epsl / (2 * erre(\n, \rn2))) * (- \bx2 * erre(\n, \rn2) + \bx2 * erre(\n, \rn2)) - (\k + 1) * \taul * \epsl * \bx2; 
                    \vy2 = Ry(\px2,\py2,\thetaangle) - (\taul * \epsl / (2 * erre(\n, \rn2))) * (- \by2 * erre(\n,\rn2) + \by2 * erre(\n, \rn2)) - (\k + 1) * \taul * \epsl * \by2; 
                    \vx3 = Rx(\px3,\py3,\thetaangle) - (\taul * \epsl / (2 * erre(\n, \rn2))) * (- \bx2 * erre(\n-1,\rn2) + \bx2 * erre(\n-1,\rn2)) - (\k + 1) * \taul * \epsl * \bx2; 
                    \vy3 = Ry(\px3,\py3,\thetaangle) - (\taul * \epsl / (2 * erre(\n, \rn2))) * (-\by2 * erre(\n-1,\rn2) + \by2 * erre(\n-1,\rn2)) - (\k + 1) * \taul * \epsl * \by2; 
                    if \k == 10 then {
                        if \n == 3 then {
                            {
                                \begin{scope}
                                    \draw[name path global = clipper, transparent] (\vcx1,\vcy1) -- (\vcx2,\vcy2) -- (\vcx3,\vcy3) -- (\vcx4,\vcy4) -- cycle;
                                    \clip (\vcx1,\vcy1) -- (\vcx2,\vcy2) -- (\vcx3,\vcy3) -- (\vcx4,\vcy4) -- cycle;
                                    \draw[name path global = last, line join = bevel] (\vx1,\vy1) -- (\vx2,\vy2) -- (\vx3,\vy3) -- cycle;
                                \end{scope}
                                \draw[intersection segments={of=clipper and last,sequence = L2}];
                            };
                        };
                    } else {
                        {\draw[line join = bevel] (\vx1,\vy1) -- (\vx2,\vy2) -- (\vx3,\vy3) -- cycle;}; 
                    };
                };
            };
        }
    \end{tikzpicture}   
    \caption{The deformed body with $b_{1} = (1/2, -\sqrt{3}/2)$ and $b_{2} = (\sqrt{3}/2, 1/2)$.} \label{figure:deformation}
\end{figure}

\begin{remark}
One of the differences between our construction and the ones in  \cite{LL} and \cite{FGS} is that the separation of the half squares introduced by the Burgers vector is $k$-dependent.
In other words, while in each square the additional horizontal shift---compared with the one above---is of one Burgers vector, the total shift in the $k$-square is the cumulative effect of all the shifts above it, since once a dislocation is added in the square above, the square below will have to accommodate it as well (see Figure~\ref{figure:deformation}).
\end{remark}

\subsection{Admissibility of $\beta$. }\label{sect:constraints}
In this section we check that the strain $\beta$ defined in \eqref{strain-beta} is admissible, namely that it satisfies assumptions \ref{regularity}--\ref{boundary-conditions}, and identify the set $S$ where the curl concentrates. We note that, from \eqref{BL:energy}, the core energy penalises the length of the interfaces between constant values of the strain that are not rank-one connected in the direction of the interface. Since we aim to prove \eqref{aim:LL}, we need to ensure that the length of such interfaces is infinitesimal. 

From our construction it is easy to see that the only lines where $\curl \beta$ may concentrate are  $\partial D_{0}$ with its translates, and $\Sigma_1 \cap \Sigma_2$. As the length of $\Sigma_1\cap \Sigma_2$ is of order one, a curl concentration there would be incompatible with the bound \eqref{aim:LL}. 
Hence we need to ensure that the values of $\beta$ across it are rank-one connected.

By imposing that the constant values of $\beta$ in $\Delta_{1,k}^{b,\bar n_1}$ and $\Delta_{2,k}^{b,\bar n_2}$ are rank-one connected along the common boundary we get the condition
$$
\left[\left(R_{-\theta}-\frac{1}{2r_{\bar n_{1}}}\, \tau\varepsilon b_1\otimes\left(\begin{matrix}\medskip
1 \\
1
\end{matrix}\right)\right)-\left(R_{\theta} -\frac{1}{2r_{\bar n_{2}}}\, \tau \varepsilon b_2\otimes\left(\begin{matrix}\medskip
1 \\
-1
\end{matrix}\right)\right)\right]\, e_2=0,
$$
which simplifies to
\begin{equation}\label{Sigma12}
\frac{1}{r_{\bar n_{1}}}\, \tau\varepsilon b_1 +\frac{1}{r_{\bar n_{2}}}\, \tau\varepsilon b_2=4\sin\theta e_1.
\end{equation}
Since $b_{1},b_{2} \in \mathbb{S}^{1}$, we may write
$$
b_1 = \left(\begin{matrix}
\cos\varphi  \\
\sin\varphi 
\end{matrix}\right), 
\quad 
b_2 = \left(\begin{matrix}
\cos\eta  \\
\sin\eta 
\end{matrix}\right)
$$
for some $\varphi, \eta \in [0, 2\pi)$.
Expanding condition \eqref{Sigma12}, we get 
$$
\begin{cases}\medskip
\displaystyle \frac{1}{r_{\bar n_{1}}} \cos\varphi + \frac{1}{r_{\bar n_{2}}} \cos\eta = \frac{4 \sin\theta}{\tau \varepsilon},\\
\displaystyle\frac{1}{r_{\bar n_{1}}} \sin\varphi + \frac{1}{r_{\bar n_{2}}} \sin\eta = 0.
\end{cases}
$$
These conditions fix the free parameters\footnote{Having fixed $r_{\bar n_{i}}$ as in \eqref{periods}, by the definition of $r_{0}$ in \eqref{squaresides}, the values of $\bar n_{i}$ are, in general, not integers. However, it can be easily ensured that $\bar n_{i}\in \mathbb{N}$ by adjusting the value of $r_{0}$ in each strip. As this procedure does not affect the results, and greatly burdens the notation, we avoid it.} $r_{\bar n_{i}}$ to be 
\begin{equation}\label{periods}
    r_{\bar n_{1}}=-\frac{\tau\varepsilon}{4}\frac{\sin(\varphi-\eta)}{\sin\theta\sin\eta},
    \qquad 
    r_{\bar n_{2}}=\frac{\tau\varepsilon}{4}\frac{\sin(\varphi-\eta)}{\sin\theta\sin\varphi}.
\end{equation}
Note that we should also ensure that $0<r_{\bar n_{i}}<L$ for $i=1,2$. While the upper bound is trivial as $\varepsilon \ll 1$, positivity is obtained when 
\begin{align}
    \frac{\sin(\varphi-\eta)}{\sin \eta} < 0,\label{condition-sign-1}\\
    \frac{\sin(\varphi-\eta)}{\sin \varphi} > 0.\label{condition-sign-2}
\end{align}
Clearly, \eqref{condition-sign-1}--\eqref{condition-sign-2} can always be satisfied up to swapping $b_{1}$ and $b_{2}$ and up to a rotation by $\pi$ of either $b_{1}$ or $b_{2}$. 

The strain $\beta$ defined in \eqref{strain-beta}, with $r_{\bar n_{i}}$ given by \eqref{periods}, satisfies \ref{regularity}, with $S = \cup_{i,k}(\partial D_{0} + t_{i,k})$, and \ref{boundary-conditions}.

We now check condition \ref{circulation}.
To compute the curl of $\beta$ on $\partial D_{0} + t_{i,k}$ it is sufficient to take $\gamma$ as a closed curve surrounding $D_{0}+t^k_i$.
With no loss of generality we can take a concentric square with side-length $2r_{n}$.
Since $\nabla I_{\Delta_{i,k}^{a,n}}$ and $\nabla I_{\Delta_{i,k}^{b,n}}$ (for $i=1,2$) are rank-one connected on their common boundary, we can take either value in the computation of the curl of $\beta$ on $\gamma$. Moreover, we can just consider the case $k=1$ since the gradients of the interpolations are independent of $k$. 
One can easily check that for $i=1$
\begin{align*}
    \int_\gamma \beta\, t ds =
    2r_{n}\left( \nabla I_{\Delta^{b,n}_1}e_2 -R_{-\theta} e_2\right)=
    -\tau\varepsilon b_1 ,
\end{align*}
and for $i=2$
\begin{align*}
    \int_\gamma \beta\, t ds &=
    2r_{n}\left(R_\theta e_2- \nabla I_{\Delta^{b,n}_2}e_2\right)= -\tau\varepsilon b_2.
\end{align*}
Hence condition \ref{circulation} is satisfied, and the strain $\beta$ is admissible. Note that our construction satisfies a stronger version of \ref{circulation}, since the circulation of $\beta$ along a closed curve $\gamma$ surrounding a single connected component of $S$ not only lies in $\tau\varepsilon \mathcal{B}$, but it is one of its generators. Heuristically, this means that each component contains a single dislocation.

\begin{remark} 
If either $\eta$ or $\varphi$ is zero, the construction presented above does not immediately work. 
Indeed, by \eqref{periods}, the spacing of the dislocations with Burgers vector $b_1$, or $b_{2}$, would be infinite.
This is due to the fact that to achieve symmetric boundary conditions only $e_1$ is needed. 
To treat this special case one has to make the simple adaptation of taking $b_{1} = b_2 = e_1$ in the construction, and consequently
\begin{equation*}
    r_{\bar n_{1}}=r_{\bar n_{2}} \coloneq \frac{\tau \varepsilon}{4 \sin\theta}.
\end{equation*}
Alternatively, one can slightly change the construction and have a single vertical strip, say $\Sigma_{1}$,  hosting the grain boundary construction.
\end{remark}

\subsection{The energy} 
We now estimate the elastic and core energy of our construction.
In doing so, we keep track of the dependence of the constants in the estimates on $\varepsilon, \theta, \tau, \lambda, \varphi, \eta$ and $L$. In particular, whenever we write simply $C$ we tacitly assume that the constant is independent of the parameters above.

\subsubsection{The elastic energy} 
For the elastic energy, note that
\begin{align*}
    \int_{\Omega}\dist^2(\beta, \SO(2))dx =\sum_{i=1}^2\sum_{k=1}^{N_i} \sum_{n=1}^{\bar n_i} \int_{(\Delta_{i,k}^{a,n}\cup \Delta_{i,k}^{b,n})\cap\Omega}\dist^2(\beta, \SO(2))dx,
\end{align*}
where we recall that $N_i = \lceil L/\ell_i \rceil$ is the number of squares $Q_i^k$ in the vertical strip $\Sigma_i$.
By the definition of $\beta$, we have 
\begin{equation}\label{estimate-dyadic}
\dist^2(\beta, \SO(2))\leq C\frac{1}{r_{n}^2} \tau^2\varepsilon^2 \quad \text{in } (\Delta_{i,k}^{a,n}\cup \Delta_{i,k}^{b,n}) \cap \Omega, \qquad  k = 1, \dots, N_{i}, \, i=1,2.
\end{equation}
Then, by \eqref{estimate-dyadic}, recalling the definition of $N_i$, we can bound
\begin{equation}\label{eq:estimate-elastic-1}
    \begin{aligned}
        \int_{\Omega}\dist^2(\beta, \SO(2)) dx &\leq C\tau^2\varepsilon^2\sum_{i=1}^2\sum_{k=1}^{N_i} \sum_{n=1}^{\bar n_i} \frac{1}{r_{n}^2}\mathcal{L}^2(D_{n}) \\
         &= C\tau^2\varepsilon^2\sum_{i=1}^2\sum_{k=1}^{N_i} \sum_{n=1}^{\bar n_i} \frac{1}{r_{n}^2}(r_{n}^2-r_{n-1}^2) \\
         &\leq C\tau^2\varepsilon^2  (\bar n_1 N_1+\bar n_2 N_2)   \\
        & \lesssim C\tau^2\varepsilon^2 L\left(  \frac{\bar n_1}{r_{\bar n_{1}}}+\frac{\bar n_2}{r_{\bar n_{2}}}\right).
    \end{aligned}
\end{equation}
In the estimate above, the last square, namely $Q_{i}^{N_{i}}$, might be cut by the boundary of $\Omega$. However, the deformation, and hence $\beta$, are in fact defined in the whole of $Q_{i}^{N_{i}}$ by interpolation (see Section \ref{strain:con}). Hence, the elastic energy of $\beta$ in $Q_{i}^{N_{i}}\cap \Omega$ can be estimated from above by the energy of its natural extension in $Q_{i}^{N_{i}}$.

We now write the right-hand side of \eqref{eq:estimate-elastic-1} more explicitly. First of all, by \eqref{form-ell}--\eqref{squaresides} and  \eqref{periods} we have
\begin{align*}
    \bar n_{1} & = \log_{2}\Big(\frac{r_{\bar n_{1}}}{r_{0}}\Big) = \log_{2} \Big(-\frac{\tau \sin(\varphi - \eta)}{4\lambda \sin\theta\sin\eta}\Big) = \log_{2}\Big(\frac{1}{\sin\theta}\Big) + \log_{2} \Big(-\frac{\tau \sin(\varphi - \eta)}{4\lambda \sin\eta}\Big),\\
    \bar n_{2} & = \log_{2}\Big(\frac{r_{\bar n_{2}}}{r_{0}}\Big) = \log_{2} \Big(\frac{1}{\sin\theta}\Big) + \log_{2} \Big(\frac{\tau \sin(\varphi - \eta)}{4\lambda \sin\varphi}\Big).
\end{align*}
By using again the expressions of $r_{\bar n_{i}}$ given in \eqref{periods}, it follows that
\begin{align*}
\frac{\bar n_1}{r_{\bar n_{1}}}+\frac{\bar n_2}{r_{\bar n_{2}}}&=
\frac{4}{\tau\varepsilon}\frac{\sin\theta}{\sin(\varphi-\eta)}\left( - \bar n_1 \sin\eta +\bar n_2\sin\varphi\right)\\
& =\frac{4}{\tau\varepsilon}\frac{\sin\theta}{\sin(\varphi-\eta)}
\log_2\left(\frac{1}{\sin\theta}\right)(\sin\varphi-\sin\eta)\\
& \hphantom{=} \, +\frac{4}{\tau\varepsilon}\frac{\sin\theta}{\sin(\varphi-\eta)}\left( - 
\sin\eta\log_2\left(-\frac{ \tau\sin(\varphi-\eta)}{4\lambda  \sin\eta}\right) + \sin\varphi\log_2\left(\frac{ \tau \sin(\varphi-\eta)}{4\lambda \sin\varphi}\right)\right).
\end{align*}
Hence, by \eqref{eq:estimate-elastic-1},
\begin{align}\label{eq:estimate-elastic}
\int_{\Omega}\dist^2(\beta, SO(2))dx &\leq C\varepsilon L \sin\theta(E_{1} - E_0\log_2(\sin\theta)),
\end{align}
where
\begin{align}
    E_{0} & \coloneq  \tau  \left( \frac{\sin \varphi - \sin \eta}{\sin(\varphi - \eta)}\right),\label{E0}\\
    E_{1} & \coloneq \tau\Big(\frac{\sin \varphi}{\sin(\varphi - \eta)} \log_{2} \Big(\frac{\tau \sin(\varphi - \eta)}{4\lambda \sin \varphi}\Big) - \frac{\sin \eta}{\sin(\varphi - \eta)} \log_{2} \Big(-\frac{\tau \sin(\varphi - \eta)}{4\lambda \sin \eta}\Big)\Big).\label{E1}
\end{align}

\subsubsection{The curl energy.} 
The conditions on the strain $\beta$ imposed in Section \ref{sect:constraints} ensure that $\curl \beta\subset S$, where $S = \cup_{i=1}^2\cup_{k=1}^{N_i}(\partial D_{0} + t_{i,k})$.
More precisely, $\curl \beta$ is concentrated on the vertical boundaries of the  translated inner squares $\partial D_{0}$. As observed above, the last square $Q_i^{N_i}$ may be cut by the boundary of $\Omega$, and therefore the curl of $\beta$ in $Q_i^{N_i}\cap \Omega$ may be concentrated in a much smaller region than $\partial D_{0}+ t_i^{N_i}$. Nevertheless, estimating the support of $\curl\beta$ with the whole of $S$ is not too crude and simplifies the computations. 
By \eqref{squaresides},
\begin{equation}\label{eq:estimate-area-core}
    \mathcal{L}^{2}(B_{\lambda \varepsilon}(\partial D_{0})) \leq C \lambda^{2}\varepsilon^{2}.
\end{equation}
Thus, by \eqref{eq:estimate-area-core} and the definition of $N_i$, the core energy in \eqref{BL:energy} can be estimated by
\begin{equation}\label{eq:estimate-core}
    \begin{aligned}
        \frac{\tau^{2}}{\lambda^{2}}\sum_{i=1}^{2}\sum_{k=1}^{N_{i}} \mathcal{L}^{2}(B_{\lambda \varepsilon}(\partial D_{0} + t_{i,k})) & \leq C\tau^{2}\varepsilon^{2}(N_{1}+N_2)\\
        & \lesssim C \tau^{2} \varepsilon^{2}L\Big(\frac{1}{r_{\bar n_{1}}} + \frac{1}{r_{\bar n_{2}}} \Big)\\
        & = C \tau \varepsilon L \sin \theta\left(\frac{\sin \varphi - \sin \eta}{\sin(\varphi - \eta)}\right)\\
        &=C\varepsilon L E_0\sin \theta,
    \end{aligned}
\end{equation}
where we used \eqref{E0}.

\subsubsection{The total energy.}
The total energy, from \eqref{eq:estimate-elastic} and \eqref{eq:estimate-core}, can be estimated from above as 
\begin{equation}\label{eq:estimate-energy}
\mathcal{F}_\e(\beta,S) \leq C \varepsilon L E_{0} \sin\theta (A - \log_{2}(\sin\theta)), 
\end{equation}
where $A:=E_1/E_0+1$, with $E_0$ and $E_1$ in \eqref{E0}--\eqref{E1}. Note that we can simplify
\begin{align}\label{def:A}
    A  = 1 + \log_{2}\Big(\frac{\tau \sin(\varphi - \eta)}{4  \lambda}\Big) + \frac{\sin \eta\log_{2}(-\sin \eta)- \sin \varphi \log_{2}(\sin\varphi) }{\sin \varphi - \sin \eta}.
\end{align}
Clearly \eqref{eq:estimate-energy} proves the upper bound \eqref{aim:LL}, for a small misorientation angle $\theta$.

\begin{remark}[The square lattice]
In the special case of a square lattice, we can fix $\varphi \in (0,\pi/2)$, and set $\eta = \varphi + 3\pi/2$.
Then the values of $r_{\bar n_{i}}$ in \eqref{periods} simplify to
$$
r_{\bar n_1} = \frac{\tau\varepsilon}{4\sin\theta\cos\varphi},
\qquad 
r_{\bar n_2} =  \frac{\tau\varepsilon}{4\sin\theta\sin\varphi}.
$$
Note that the spacings between dislocations, $2r_{\bar n_i}$, match perfectly the ones computed in \cite[Eqs.\ (3)--(4)]{RS} by Read and Shockley (we recall that in \cite{RS} the angle of misorientation is $\theta / 2$). As for the energy estimate \eqref{eq:estimate-energy}, the constants $E_0$ and $A$ (defined in \eqref{E0} and \eqref{def:A}) 
reduce to 
\begin{align}
E_0 &= \tau \big(\sin\varphi+\cos\varphi\big),\label{E0-square}\\
A &= 1 + \log_2\left(\frac{\tau}{4\lambda}\right) -\frac{\cos\varphi\log_2(\cos\varphi)+ \sin\varphi\log_2(\sin\varphi)}{\sin\varphi+\cos\varphi}.\label{A}
\end{align}
The constant $E_0$ in \eqref{E0-square} is in complete agreement with the formula derived by Read and Shockley\footnote{Recall that in \cite{RS} the angle of misorientation is $\theta / 2$. However, up to changing the constant $C$ in \eqref{eq:estimate-energy}, no changes to $E_{0}$ or $A$ are needed.}, see \cite[Eq.\ (6)]{RS}.
Hence, at leading order, the energy of our construction gives the optimal scaling with the correct constant. As for the constant $A$ in the lower order term of the energy, the corresponding one in \cite[Eqs.\ (7)--(8)]{RS} is 
\begin{equation}\label{ARS}
    A_{\text{RS}} = 1 + \log_2\left(\frac{\tau}{2\pi\lambda}\right) - \frac{\sin(2\varphi)}{2} -\frac{\cos\varphi\log_2(\cos\varphi)+ \sin\varphi\log_2(\sin\varphi)}{\sin\varphi+\cos\varphi}.
\end{equation}
The presence of of $2\pi$ in \eqref{ARS} instead of $4$ in \eqref{A} can be easily explained by the fact that the core region in our construction is modelled  as a square instead of a ball.
The term $- \sin(2\varphi)/2$ in \eqref{ARS} is missing in \eqref{A} since our modelling of the core energy is too crude. In particular, the core energy in \eqref{BL:energy} does not depend on the Burgers direction and hence on the lattice orientation. Note that for $\varphi=0$ the constants are in complete agreement.
\end{remark}

\bigskip

    \noindent
    \textbf{Acknowledgements.} 
    LS is indebted to David Bourne and Myrto Galanopoulou for the useful discussions on the Read and Shockley paper, and acknowledges support by the EPSRC under the grants EP/V00204X/1 and EP/V008897/1. 
    EGT is funded by the Deutsche Forschungsgemeinschaft (DFG, German Research Foundation) under Germany's Excellence Strategy EXC 2044/2 –390685587, Mathematics Münster: Dynamics–Geometry–Structure.

\end{document}